\documentclass[12pt]{article}
\usepackage{amssymb}
\usepackage{amsmath}
\usepackage{amsbsy}
\usepackage{amsthm}
\usepackage{amscd}

\newcommand{\pt}{\partial}

\newcommand{\Pf}{\noindent {\em Proof.} }

\newtheorem{theo}{Theorem}[section]
\newtheorem{rem}[theo]{Remark}
\newtheorem{prop}[theo]{Proposition}
\newtheorem{lem}[theo]{Lemma}
\newtheorem{notat}[theo]{Notation}
\newtheorem{term}[theo]{Terminology}
\newtheorem{cor}[theo]{Corollary}
\newtheorem{defi}[theo]{Definition}
\newtheorem{defi-lem}[theo]{Definition-Lemma}

\newtheorem{ex}[theo]{Example}
\newtheorem{exs}[theo]{Examples}

\newtheorem{ass}[theo]{Assumption}

\def \C{{\bf C}}
\def \R{{\bf R}}
\def \Z{{\bf Z}}

\def \D{{\cal D}}

\def \E{{\cal E}}
\def \F{{\cal F}}

\def \B{{\cal B}}
\def \cC{{\cal C}}

\def \L{{\cal L}}
\def \K{{\cal K}}
\def \M{{\cal M}}

\def \Q{{\cal Q}}
\def \R{{\bf R}}
\def \cR{{\cal R}}

\def \O{{\cal O}}
\def \T{{\cal T}}
\def \V{{\cal V}}

\def \Cst{{$C^*$}}

\def \ind{{\rm ind}}

\def \mod{{\rm mod\,}}

\def\im{{\rm Im}\,}

\def \ext{{\rm ext\,}}
\def \inter{{\rm \,int\,}}
\def\supp{{\rm supp\,}}

\def\hotimes{\hat\otimes}
\def\Ct{Cl_\tau}

\def\CG{Cl_{\Gamma}}

\def\CtG{Cl_{\tau\oplus\Gamma}}

\def\c*G{C^*_\Gamma}
\def\Slf{\gtS_{lf}}
\def\Str{\gtS_{tr}}
\def\Ind{{\rm Ind}}
\def\Hom{{\rm Hom}}
\def\End{{\rm End}}

\def\l2lf{{L^2_{lf}}}

\def\Kl{{\cal K}_{loc}}

\def\gtb{\mathfrak{b}}

\def\gtc{\mathfrak{c}}
\def\gtd{\mathfrak{d}}
\def\gtM{\mathfrak{M}}
\def\gtm{\mathfrak{m}}

\def\gtS{\mathfrak{S}}
\def\gtg{\mathfrak{g}}
\def\gtG{\mathfrak{G}}
\def\gtD{\mathfrak{D}}
\def\gtd{\mathfrak{d}}
\def\gtf{\mathfrak{f}}

\def\gtG{\mathfrak{G}}

\def\cqfd{{\vbox{\hrule height 5pt width 5pt}\quad}\penalty -10}
\def\kkprod{\times\hspace{-15pt}\times}

\begin{document}

\title{$K$-theory and index theory\\
on manifolds with a proper Lie group action\\
{\small (Preliminary version 3; \; May 2024)}}

\date{}

\author{Gennadi Kasparov}

\maketitle

\begin{abstract} The paper is devoted to the index theory of orbital and transverse elliptic operators on manifolds with a proper Lie group action. It corrects errors of my previous paper on transverse operators and contains new results. Throughout the paper, we use the operator integration method in constructing pseudo-differential operators. 
\end{abstract}

\section{Introduction}

This is a continuation of my previous work \cite{Ka16}. In the present paper, additional $K$-theoretic tools are developed in order to treat index theory for leaf-wise (i.e. orbital) and transverse elliptic operators on manifolds with a proper Lie group action. Under a group action we understand a smooth, proper, isometric Lie group action. The two index theories, orbital and transverse, are very much intertwined and interdependent, and will be treated together. The proofs of index theorems for elliptic, t-elliptic and leaf-wise elliptic operators (based on the $KK$-theoretic approach) are given in sections 9, 10, 11 of the present paper.

One of the new features of the present paper is a coarse approach to the pseudo-differential calculus. A significant part of index theory uses only the PDO calculus \emph{modulo compact operators}. The coarse PDO calculus that we propose allows to treat those index theory problems which do not require very elaborate analytical tools. This is not a replacement of H\" ormander's PDO calculus. It goes along with the H\" ormander's calculus and simplifies it in many cases. 

The main index theory results for the transversally elliptic case were presented in sections 6, 7 and 8 of \cite{Ka16}. Unfortunately, there were errors in sections 7 and 8 of \cite{Ka16}. The corrections are given in the present paper. The main index theorem 8.18 of \cite{Ka16} is correct and is reproved in section 10 below.

The theory of leaf-wise operators (when the leaves are orbits of a Lie group action) is developed in the present paper from the basic definitions to the final index theorem. Although there are some common points in the group action case and the general case of leaf-wise operators on singular foliations (\cite{A-Sk II, A-Sk III}), the results in the group action case are more concrete.  

\bigskip

\noindent{\bf\large Brief contents:}

\medskip

\noindent-- The basics of the coarse PDO calculus: sections 2--3.

\smallskip

\noindent-- \Cst-algebras associated with a group action: section 4.

\smallskip

\noindent-- Leaf-wise symbols and the construction of leaf-wise operators: section 5.

\smallskip

\noindent-- Poincare duality and $K$-theory of symbol algebras: sections 6--7.

\smallskip

\noindent-- Elliptic symbols and index groups: section 8.

\smallskip

\noindent-- Index theorems for elliptic, t-elliptic and leaf-wise elliptic operators: sections 9, 10, 11.

\smallskip

\noindent-- Examples: t-elliptic and leaf-wise operators: section 12.

\smallskip

\noindent-- Comments on my article \cite{Ka16}: Appendix.

\bigskip

\noindent{\bf\large Notational references:}

\medskip

We will use certain definitions from \cite{Ka88,Ka16}. Some (but not all) of these definitions will be repeated in the present paper. For  the reader's convenience, we give here the list of references for these definitions. 

\medskip

\noindent-- $C_0(X)$-algebras and $\cR KK(X;A,B)$ groups - \cite{Ka88}, sections 1, 2.

\smallskip

\noindent-- The algebra of $G$-invariant elements $B^G$ for a proper $G$-algebra $B$ - \cite{Ka88}, 3.2.

\smallskip

\noindent-- Clifford algebras of vector bundles, especially $\Ct(X)$ for a manifold $X$, and the Dirac element $[d_X]\in K^0(\Ct(X))$ - \cite{Ka88}, section 4, or \cite{Ka16}, section 2. 

\smallskip

\noindent-- Basic definitions for transversally elliptic operators - \cite{Ka16}, section 6.

\smallskip

\noindent-- Bott and Dirac operators on tangent spaces of a manifold $X$; elements $[\B_\xi], [\B_{\xi,\Gamma}], [d_\xi]$ - \cite{Ka16}, 2.5, 2.6, 7.2.

\smallskip

\noindent-- Dirac element $[d_{X,\Gamma}]\in K^0(C^*(G,\CtG(X)))$ - \cite{Ka16}, 8.8.

\smallskip

\noindent-- Dolbeault elements $[\D_X], [\D_{X,\Gamma}]$ and $[\D^{cl}_{X,\Gamma}]$ for the tangent manifold $TX$ - \cite{Ka16}, 2.8 and 8.17.

\smallskip

\noindent-- Clifford symbol for elliptic and transversally elliptic operators - \cite{Ka16}, 3.8 and 8.11.

\smallskip

\noindent-- Tangent Clifford symbol for transversally elliptic operators - \cite{Ka16}, 8.13.

\section{Operator integration}

We will use operator integration for the construction of pseudo-differential operators. We will need two types of integration described in the two following subsections.

\subsection{Riemann integration}

This kind of operator integration is essentially described in section 3 of \cite{Ka75}. We will add here a few more details. 

Let $X$ be a second countable, locally compact, $\sigma$-compact space, $D$ a \Cst-algebra, and $\phi:C_0(X)\to \M(D)$ a homomorphism such that $\phi(C_0(X))\cdot D$ is dense in $D$. Then $\phi$ extends to a unital homomorphism $\phi:C_b(X)\to \M(D)$, where $C_b(X)$ is the \Cst-algebra of all bounded continuous functions on $X$. We will denote the set of all compactly supported continuous functions on $X$ by $C_c(X)$.

\begin{defi} The support of an element $F\in \M(D)$ is the smallest closed subset of $X\times X$, denoted $\supp(F)$, such that for any $a\in C_c(X)$, one has: $\phi(a)F=0$ as soon as $(\supp(a)\times X)\cap\supp(F)=\emptyset$, and $F\phi(a)=0$ as soon as $(X\times\supp(a))\cap\supp(F)=\emptyset$. If $\supp(F)$ is compact, $F$ will be called compactly supported. An element $F\in \M(D)$ will be called properly supported if both projections $p_1:\supp(F)\to X$ and $p_2:\supp(F)\to X$ are proper maps. 
 \end{defi}
 
 \begin{rem} {\rm Suppose $F$ is properly supported and $a\in C_c(X)$. Then both $\phi(a)F$ and $F\phi(a)$ are compactly supported. Indeed, it is easy to check that $\supp(\phi(a)F)\subset p_1^{-1}(\supp(a))\cap \supp(F)$, which is compact. Similarly for $F\phi(a)$.}
\end{rem}

\begin{defi} We will call an element $F\in \M(D)$ locally compact if for any $a\in C_0(X)$, both $\phi(a)F$ and $F\phi(a)$ belong to $D$. The set of locally compact elements will be denoted $D_{loc}$. In the case of $D=\K(H)$ for a Hilbert space $H$, the notation for locally compact elements will be $\Kl(H)$.

We denote by $Q_{C_0(X)}(D)$ the subalgebra of $\M(D)$ consisting of elements $T\in \M(D)$ which commute with $\phi(C_0(X))$ modulo $D$. The algebra $D_{loc}$ is a two-sided ideal in $Q_{C_0(X)}(D)$.
\end{defi}

Given a bounded norm-continuous map $F:X\to \M(D)$ such that $F(x)$ commutes with $\phi(C_0(X))$ modulo $D$, we will construct a Riemann type operator integral $\int_X F(x)d\phi\in Q_{C_0(X)}(D)/D_{loc}$ which has the following properties (cf. \cite{Ka75}, section 3, theorem 1):

\begin{theo} $1^{\rm o}$. If $||F(x)|| \le c$ for all $x\in X$ then $\int_X F(x)d\phi \le c$. 

$2^{\rm o}$. The integral is additive, multiplicative, and 
$\int_X F^*(x)d\phi=(\int_X F(x)d\phi)^*$.  

$3^{\rm o}$. If $F$ is a scalar function $F(x)=f(x)\cdot 1$, where $f\in C_b(X)$, then $\int_X F(x)d\phi=\phi(f)$.

$4^{\rm o}$. The integral is functorial in $X$: if there is a proper continuous map $h:Y\to X$, and $\tilde F=F\cdot h:Y\to\M(D)$, $\phi=\psi \cdot h^*: C_0(X)\to C_0(Y)\to \M(D)$, then $\int_X F(x)d\phi=\int_Y \tilde F(y) d\psi$.

$5^{\rm o}$. $\int_X F(x)d\phi$ can be lifted to $Q_{C_0(X)}(D)$ as a properly supported element.
\end{theo}

The following simple lemma (\cite{Ka75}, section 3, lemma 1) provides the necessary estimates for the proof:

 \begin{lem} Let $B$ be a \Cst-algebra, and elements $F_1,...,F_n;\alpha_1,...,\alpha_n\in B$ satisfy the following assumptions: $\sum_{i=1}^n \alpha_i^*\alpha_i=1$, and $||F_i||\le c$ for any $i$. Then $||\sum_{i=1}^n \alpha_i^*F_i\alpha_i||\le c$. 
 \end{lem}
 
 \Pf Using a faithful representation of $B$ in a Hilbert space $H$, we get for any $\xi,\eta\in H$: 
 $$|(\sum_i \alpha_i^*F_i\alpha_i(\xi),\eta)|=|\sum_i (F_i\alpha_i(\xi),\alpha_i(\eta))| \le \sum_i ||F_i|| \cdot ||\alpha_i(\xi)|| \cdot ||\alpha_i(\eta)||$$ 
 $$\le c (\sum_i||\alpha_i(\xi)||^2)^{1/2}(\sum_i||\alpha_i(\eta)||^2)^{1/2}=c||\xi|| ||\eta||. \;\cqfd $$

 \noindent\emph{Proof of the theorem.} The basic idea of the integral is the following. Let us assume that $X$ is compact. Let $B$ be the \Cst-subalgebra of $\M(D)/D$ generated by $\phi(C(X))$ and all elements $F(x),\, x\in X$. Then there is a unital homomorphism: $C(X)\to B$ which maps $C(X)$ to the center of $B$. The map $F:X\to B$ represents an element of $C(X,B)\simeq C(X)\otimes B$. Consider the multiplication homomorphism: $C(X)\otimes B\to B$. The image of $F$ under this homomorphism is by definition $\int_X F(x)d\phi\in B\subset \M(D)/D$. 
 
 Now let us translate this into the Riemann operator integration context. We will still continue to assume for the moment that $X$ is compact. The integral is constructed in the following way. Taking a finite covering $\{U_i\}$ of $X$, points $x_i\in U_i$, and a partition of unity $\sum_i \alpha_i^2(x)=1$ associated with $\{U_i\}$ (all functions $\alpha_i$ are non-negative), we consider the integral sum: $\Sigma(\{U_i\},\{\alpha_i\},\{x_i\})=\sum_i \phi(\alpha_i)F(x_i)\phi(\alpha_i)$. We assume that $||F(x)||\le c$ for all $x\in X$, so by lemma 2.5 the norm of the integral sum is $\le c$. 
 
 We will call $\{U_i\}$ an $\epsilon$-covering if for any $i$ and any $x,y\in U_i$, one has $||F(x)-F(y)||\le\epsilon$. To verify that the integral sums for two $\epsilon$-coverings $\{U_i\}$ and $\{V_j\}$ differ in norm no more than by $2\epsilon$ in $\M(D)/D$, we form the covering $W_{i,j}=U_i\cap V_j$. Let $\sum_j\beta_j^2=1$ be the partition of unity for $\{V_j\}$, and set $\gamma_{i,j}=\alpha_i\beta_j,\, x_{i,j}\in W_{i,j}$. Then we have the following estimate in $\M(D)/D$:
 $$||\sum_{i,j}\phi(\gamma_{i,j})(F(x_i)-F(x_{i,j}))\phi(\gamma_{i,j})||\le \epsilon$$
by lemma 2.5. Here 
$$\sum_{i,j}\phi(\gamma_{i,j})F(x_i)\phi(\gamma_{i,j})=\sum_i\phi(\alpha_i) F(x_i)\phi(\alpha_i)$$ 
modulo $D$ because $\sum_j \beta_j^2=1$. 
The integral is the limit in $\M(D)/D$ of the integral sums for all $\epsilon$-coverings when $\epsilon\to 0$. For a compact $X$, this gives the existence.  

In the general case of a non-compact $X$, let $\{U_i\}$ be a locally finite covering of $X$ and $\sum_i \alpha_i^2(x)=1$ the corresponding partition of unity. Then, by lemma 2.5, the norms of operators $\F_m=\sum_{i=1}^m \phi(\alpha_i) F(x_i) \phi(\alpha_i)$ are uniformly bounded in $m$. We claim that the sums $\F_m$ converge strictly in $\M(D)$ when $m\to \infty$. Indeed, for any bounded approximate unit $\{u_k\}\subset C_0(X)$ consisting of functions with compact support in $X$, our initial assumption on $\phi$ says that $\{\phi(u_k)\}$ converges strictly in $\M(D)$. Since $\F_m$ is uniformly bounded in $m$, it is enough to show that both sums $\F_m\phi(u_k)$ and $\F^*_m\phi(u_k)$ converge for any fixed $k$ when $m\to \infty$. This is true because all $u_k$ have compact support, so $\alpha_iu_k=0$ for large $i$.

Now we can take any locally finite $\epsilon$-covering of $X$ and form the integral sum $\F=\sum_{i=1}^\infty \phi(\alpha_i) F(x_i) \phi(\alpha_i)$. Modulo $D_{loc}$ any two such sums for two different $\epsilon$-coverings will differ in norm by $\le 2\epsilon$ (because all  corresponding finite sums $\F_m$ differ by $\le 2\epsilon$, as shown above). So in $Q_{C_0(X)}(D)/D_{loc}$ these sums converge in norm when $\epsilon\to 0$. The lifting of the limit to $Q_{C_0(X)}(D)$ is denoted $\int_X F(x)d\phi$. 

The proof of all properties listed in the theorem (except $5^{\rm o}$) for the case of a non-compact $X$ is the same as for a compact $X$ (see \cite{Ka75}, section 3, theorem 1). For example, for the multiplicativity property, one needs to estimate the difference modulo $D_{loc}$ of the two sums:
$$\sum_{i,j}\phi(\alpha_i)F_1(x_i)\phi(\alpha_i)\phi(\alpha_j)F_2(x_j)\phi(\alpha_j)-\sum_{i,j}\phi(\alpha_i)F_1(x_i)\phi(\alpha_i)\phi(\alpha_j)F_2(x_i)\phi(\alpha_j).$$
Both sums converge in the strict topology, so one needs only to show that all finite portions of these sums differ modulo $D_{loc}$ by less than $2\epsilon c$ if $\{U_i\}$ is an $\epsilon$-covering for $F_2$ and $||F_1(x)||\le c$ for all $x\in X$. Obviously, we can leave only those summands for which $U_i\cap U_j\ne\emptyset$. Then modulo $D_{loc}$, $||F_2(x_i)-F_2(x_j)||\le 2\epsilon$, and the estimate of the difference of those two sums $\sum_{1\le i,j\le n}$ modulo $D$ comes from lemma 2.5. 

Concerning property $5^{\rm o}$, we can assume that all our coverings $\{U_i\}$ consist of open sets with compact closure. This means that the operator $S=\sum_i \phi(\alpha_i) F(x_i) \phi(\alpha_i)$ is properly supported. Indeed, it is easy to check that $\supp(S)\subset \cup_i (\supp(\phi(\alpha_i))\times \supp(\phi(\alpha_i)))$ in $X\times X$. Here $\supp(\phi(\alpha_i))$ means the following: The homomorphism $\phi$ maps $C_0(X)$ onto a commutative \Cst-subalgebra in $\M(D)$. The spectrum of this commutative \Cst-algebra can be identified with a closed subset $Y\subset X$. Then $\supp(\phi(\alpha_i))$ is a compact subset in $Y$, and hence in $X$.

It is clear from the construction of the integral that once we have chosen one covering $\{U_i\}$, all subsequent coverings can be chosen as the intersections of this one with other coverings like $\{V_j\}$ above. Also the corresponding partitions of unity can be chosen as the products $\{\alpha_i\beta_j\}$. So if we fix $\{\alpha_i\}$ once and for all and change only $\{V_j,\beta_j\}$, then the resulting integral will have the form: $\sum_i \phi(\alpha_i)(\int_{U_i} F(x)d\phi)\phi(\alpha_i)$ which obviously lifts to a properly supported element. \cqfd
 
 \begin{cor} In the assumptions of the theorem, suppose that $||F(x)||\to 0$ when $x\to\infty$ in $X$. Then the sums $\F_m$ in the above proof of the theorem converge uniformly in $\M(D)$.
\end{cor}
  
 \Pf For any $\epsilon>0$, we can choose $m$ and $n$ large enough so that for any $x$ in the supports of all functions $\alpha_i$, for $i$ in the interval $[m,n]$, we have $||F(x)||\le\epsilon$. Then $||\F_m-\F_n||\le \epsilon$ by lemma 2.5. \cqfd

  \subsection{Group averaging}
  
  We keep all assumptions of the previous subsection, but assume in addition that a locally compact, second countable group $G$ acts on $X$ properly, $D$ is a $G$-algebra, and the homomorphism $\phi:C_0(X)\to \M(D)$ is $G$-equivariant. In integration over $G$, we will use the left Haar measure.
 
 \begin{lem} If $F:X\to \M(D)/D$ is $G$-equivariant, then $I(F)=\int_X F(x) d\phi$ is $G$-invariant modulo $D$ and $G$-continuous in norm (i.e. the map $G\to \M(D): g\mapsto g(I(F))$ is norm-continuous). 
 \end{lem}
 
 \Pf Because the $G$-action transforms an $\epsilon$-covering into another $\epsilon$-covering, it is clear that $I(F)$ is $G$-invariant modulo $D$.  The last assertion follows from \cite{Thom_prepr}, 1.1.4. \cqfd
 
In order to make $I(F)$ exactly $G$-invariant we need averaging over $G$. Essentially, we will use the averaging method of \cite{CM:L^2}, proposition 1.4. We will adapt it to the generality that we need.
 
\begin{prop} Let $T\in \M(D)$ be an operator with support in a set $L\times L$, where $L$ is a compact subset of $X$. Then one can define the average of $T$ over $G$, denoted $Av_G(T)$ or $\int_G g(T)dg$, as the limit of integrals $\int_C g(T)dg$, in the strict topology of $\M(D)$, over the increasing net of all compact subsets $C$ of $G$, where $dg$ denotes the Haar measure of $G$. Moreover, $||Av_G(T)||\le c ||T||$, where $c$ depends only on $L$.
 \end{prop}
 
 For the proof we need the following lemma:
 
 \begin{lem} If an operator $T\in \M(D)$ has the property that $T^*g(T)=0$ and $Tg(T^*)=0$ for any $g$ outside of a compact set $K\subset G$, then for any compact subset $C\subset G$, $(\int_C g(T)dg)^*(\int_C g(T)dg)\le |K|^2||T||^2$, where $|K|$ is the Haar measure of $K$.
 \end{lem}
 
 \Pf Consider $[(\int_C g(T)dg)^*(\int_C g(T)dg)]^n$ as the $2n$-fold integral 
 $$\int_C...\int_C g_1(T^*)g_2(T)...g_{2n-1}(T^*)g_{2n}(T)dg_1...dg_{2n}.$$
 Because of our assumption on $T$, this multiple integral actually goes over the subset of $C\times...\times C$ such that $g_i^{-1}g_{i+1}\in K$ for all $i$. So we can rewrite it as a repeated integral:
 $$\int_C\int_K...\int_K g_1(T^*h_2(Th_3(T^*h_4(T...h_{2n-1}(T^*h_{2n}(T)...)dg_1dh_2...dh_{2n},$$
where $h_{i+1}=g_i^{-1}g_{i+1}$. The latter integral is estimated as $\le |C|\cdot |K|^{2n-1}||T||^{2n}$, therefore $[(\int_C g(T)dg)^*(\int_C g(T)dg)]^n\le |C|\cdot |K|^{2n-1}||T||^{2n}$. Taking the $n$-th root of both sides and letting $n\to\infty$, we get the result. \cqfd

\medskip

\noindent\emph{Proof of the proposition.} If $a\in C_c(X)$ is a function equal to $1$ on $L$, then $T\phi(a)=T$ and $T^*\phi(a)=T^*$. Because the $G$-action on $X$ is proper, there is a compact set $K\subset G$ such that $\supp(a)$ does not intersect with any $g(L)$ for any $g\in G, \; g\notin K$. So for such $g$ we have: $T^*g(T)=T^*\phi(a)g(T)=0$, and similarly for $Tg(T^*)$. In view of the previous lemma, $||\int_C g(T)dg||\le |K| \cdot ||T||$ for any compact subset $C\subset G$, and the same for $T^*$. 

Let us denote $\int_C g(T)dg$ by $I(C)$. We need to show that for any $d\in D$, the integrals $d\cdot I(C)$ and $d\cdot I(C)^*$ converge over the net of all compact $C$. Clearly, this is true if we replace $d\in D$ with $\phi(a)$, where $a\in C_c(X)$, because if we define $C_{max}(a,T)$ as the maximal compact subset of $G$ such that $\supp(a)\cap g(L)\ne \emptyset$, then both integrals will not depend on $C$ as soon as $C_{max}(a,T)\subset C$. 

Given $d\in D$ and $\epsilon>0$, we can find $a\in C_c(X)$ such that $||d-d\cdot\phi(a)||\le \epsilon$ and $||d-\phi(a)\cdot d||\le \epsilon$. By the previous lemma, both $I(C)$ and $I(C)^*$ are bounded as functions of $C$, more precisely, $||d\cdot I(C)-\phi(a)\cdot d\cdot I(C)||\le \epsilon |K|\cdot ||T||$, and similarly for $I(C)^*$. This means that $d\cdot I(C)$ varies with $C$ no more than by $\epsilon |K|\cdot ||T||$ when $C_{max}(a,T)\subset C$. This proves the convergence. The norm estimate is also clear. \cqfd

 \section{Coarse PDO calculus}

Index theory deals with what is usually called `elliptic' operators. These are some kind of bounded (`zero order') operators, invertible modulo some kind of `negative order' operators. In the simplest case, one has a compact manifold $X$ and a representation of $C(X)$ in a Hilbert space $H=L^2(X)$. `Zero order' operators are bounded operators on $H$ which commute with the action of $C(X)$ modulo $\K(H)$, `negative order' operators are just compact operators. If $X$ is locally compact, the requirement for `zero order' operators is the same: commutation with the action of $C_0(X)$ modulo $\K(H)$. `Negative order' operators $T$ are those which satisfy the conditions: $T\cdot C_0(X)\subset \K(H)$ and $C_0(X)\cdot T\subset \K(H)$. (So `negative order' operators are exactly what was called `locally compact' operators in definition 2.3.) Index theory seeks invariants of operators of `zero order' modulo `negative order'. There is an obvious similarity with the basic $KK$-theory (cf. \cite{Ka88}, definition 2.1).

We will use this scheme (`zero order' modulo `negative order' operators) in our coarse pseudo-differential (PDO) calculus. This calculus will correspond to the H\" ormander $\rho=1,\delta=0$ PDO calculus. In fact, H\" ormander discussed such kind of `simple' PDO calculus in \cite{Hor} at the end of subsection 2.1, but based on his definition of PDOs. Our approach is different and aimes at simplifying as much as possible the definition of a PDO for zero order operators. 

There are positive and negative sides of this approach. On the positive side, the whole PDO theory gets much simpler. On the negative side, this approach is not very convenient in treating differential operators and their parametrices. But when necessary, one can always shift back to H\" ormander's approach which is perfectly compatible with the coarse one that we propose.

Our coarse PDO theory uses the operator integration of section 2. In all constructions, the Riemannian metric of the manifold and a Hermitian metric on vector bundles will play an important role. 

In what follows, $X$ will be a complete Riemannian manifold. Our setting may include a proper isometric action of a locally compact group $G$ on $X$. We denote by $T(X)$ the tangent bundle of $X$, by $T^*(X)$ the cotangent bundle and by $p:T^*(X)\to X$ the projection. We will often identify $T(X)$ and $T^*(X)$ via the Riemannian metric of $X$. The tangent manifold will be denoted $TX$, and the projection $TX\to X$ also by $p$.

\subsection{Preliminaries}

We start with the following theorem (which most likely exists somewhere in the literature).

\begin{theo} Let $f=f(\xi)$ be a bounded, differentiable function on $\R^n$ such that all its first derivatives vanish at infinity. Denote by $\Phi$ the operator of Fourier transform on $L^2(\R^n)$, by $f$ the operator of multiplication by the function $f$ on $L^2(\R^n)$, and by $F$ the operator $\Phi^{-1}f\Phi$. Then for any function $a=a(x)\in C_0(\R^n)$ considered as a multiplication operator on $L^2(\R^n)$, the commutator $[F,a]$ is a compact operator on $L^2(\R^n)$.

More generally, the assertion remains true in the situation with a compact parameter space $Z$. More precisely, if $L^2(\R^n)$ is replaced with $L^2(\R^n)\otimes C(Z)$ and $f(\xi)$ with $f(\xi,z)\in C_b(\R^n\times Z)$, continuous in $z$ uniformly in $\xi$, satisfying the same assumption on its first derivatives in $\xi$ (uniformly in $z\in Z$), then for any $a\in C_0(\R^n\times Z)$, the commutator $[F,a]$ belongs to $\K(L^2(\R^n))\otimes C(Z)$.
\end{theo}

\Pf We will prove the Fourier-dual assertion. Recall that $\Phi^{-1}C_0(\R^n)\Phi=C^*(\R^n)$, the \Cst-algebra of the abelian group $\R^n$. The algebra $C^*(\R^n)$ contains the dense subalgebra $C_c(\R^n)$ (compactly supported continuous functions) with convolution as multiplication. We need to prove that for any $b\in C_c(\R^n)$, the commutator $[b,f]\in \K(L^2(\R^n))$.

The assumption on the first derivatives of $f$ implies that if $x\to\infty$ (or $y\to\infty$) and $||x-y||$ remains bounded, then $|f(y)-f(x)|\le \int_0^1 |\partial/\partial t f(x+t(y-x))|dt\to 0$. The commutator $[b,f]$ is an integral operator with the kernel $k(x,y)=b(x-y)(f(x)-f(y))$. Since $b$ has compact support, we obviously get $\int |k(x,y)|dx\to 0$ when $y\to\infty$ and $\int |k(x,y)|dy\to 0$ when $x\to\infty$. The Schur lemma \cite{Hor_book}, 18.1.12, easily implies that the integral operator with the kernel $k$ is compact. 

The proof of the generalized version of the statement (with the parameter space $Z$) is the same. In fact, it is enough to work with $f(\xi,z)=f_1(\xi)f_2(z)$ and $a(x,z)=a_1(x)a_2(z)$. \cqfd

\begin{prop} $1^{\rm o}$ The assertion of theorem 3.1 remains true if $f$ is bounded and measurable, but differentiable only outside of a compact subset of $\R^n$, with all first derivatives of $f$ vanishing at infinity. Moreover, the norm of the operator $F\,\mod \K(L^2(\R^n))$ does not exceed $\limsup_{\xi\to\infty}|f(\xi)|$.

$2^{\rm o}$ Keeping all assumptions of theorem 3.1 concerning $f$ and $F$, let also $a\in C_0(\R^n-\{0\})$. Then the product of operators $F\cdot a$ belongs to $\K(L^2(\R^n))$. In particular, $F\cdot(a-b)$ is a compact operator if $a,b\in C_0(\R^n)$ and $a(0)=b(0)$. 
\end{prop}

\Pf $1^{\rm o}$ We can write $f=f_0+f_1$ where $f_0$ is bounded,  measurable and has compact support, and $f_1$ satisfies the assumptions of the theorem. For any $a\in C_0(\R^n)$, the products $\Phi^{-1}f_0\Phi\cdot a$ and $a\cdot\Phi^{-1}f_0\Phi$ belong to $\K(L^2(\R^n))$ by the Rellich lemma. So the assertion of the theorem remains true. For the last statement, we can take $f_0$ with as large compact support as we want. The norm of $F\,\mod \K$ depends only on $\sup |f_1|$. 

$2^{\rm o}$ Again we will prove the Fourier-dual assertion: $f\cdot \hat a\in \K(L^2(\R^n))$, where $\hat a=\Phi^{-1}a\Phi$ is the Fourier-dual operator. By the Stone theorem, the algebra $C_0(\R^n-\{0\})$ is  generated by functions $x_k(\exp(-||x||^2/2))$ with $1\le k\le n$. The Fourier-dual operators for these functions are (up to a scalar multiple) the same functions (in $\xi$ instead of $x$), which can also be written as $\partial/\partial \xi_k\exp(-||\xi||^2/2)$. If $\hat a$ is written in this form, then $f\hat a=f\cdot\pt/\pt\xi_k(\exp(-||\xi||^2/2))=\pt/\pt\xi_k(f)\cdot\exp(-||\xi||^2/2)\in \K(L^2(\R^n))$. \cqfd

\begin{prop} Keeping the assumptions of theorem 3.1 concerning the function $f$ and the operator $F$, let $V$ be a finite-dimensional complex Hermitian vector space, $h$ and $k$ - continuous maps $\R^n\to \L(V)$, such that $h=h(x)$ is linear in $x$ and $(h(x))^*=h(x)$ for any $x$, $k=k(x)$ is bounded, and both $||(h^2+1)^{-1}k||$ and $||(h^2+1)^{-1}hk||$ vanish at infinity in $x$. Then we have: $[F,h(h^2+1)^{-1/2}]k\in \K(L^2(\R^n)\otimes V)$. 
\end{prop}

\Pf We will use the integral presentation: 
$$h(h^2+1)^{-1/2}=2/\pi \int_0^\infty h(1+h^2+\lambda^2)^{-1}d\lambda.$$
The commutator $[F,h(h^2+1)^{-1/2}]$ equals 
$$2/\pi \int_0^\infty (1+h^2+\lambda^2)^{-1}\cdot ((1+\lambda^2)[F,h]+h[F,h]h)\cdot(1+h^2+\lambda^2)^{-1}d\lambda.$$
Since $h$ is a linear function, the commutator $[F,h]$ belongs to $C^*(\R^n)\otimes \L(V)$. (Indeed, the Fourier-dual assertion says that the commutator of $f$ with the first order differential operator $\hat h$ equals to the multiplication by a function vanishing at infinity.) Therefore the latter integral expression converges in norm, and because $(h^2+\lambda^2+1)^{-1}k$ and $(h^2+\lambda^2+1)^{-1}hk$ both belong to $C_0(\R^n)\otimes\L(V)$, the assertion follows from the Rellich lemma. \cqfd

\subsection{Symbols}

\begin{lem} Let $\sigma(x,\xi)$ be a function on $K\times \R^n$, where $K$ is a compact subset of $\R^n$, and let $\{\psi_x\}:(x,\xi)\mapsto (x,\psi_x(\xi))$ be a norm-continuous family of invertible linear maps of $K\times \R^n$ into itself.

$1^{\rm o}$ Assume that $\sigma(x,\xi)$ is continuous in $x$ uniformly in $\xi$. 

$2^{\rm o}$ Also assume that $\sigma(x,\xi)$ is differentiable in $\xi$, and for the exterior derivative $d_\xi$, there is an estimate: $||d_\xi \sigma(x,\xi)||\le C\cdot (1+||\xi||)^{-1}$ with the constant $C$ which does not depend on $(x,\xi)$. 

Then $\sigma(x,\psi_x(\xi))$ satisfies the same two conditions as $\sigma(x,\xi)$ (with a different constant $C$). Moreover, if the assumptions concerning $\sigma$ hold only outside of a compact subset in $\xi\in\R^n$, the assertion remains true outside of a (possibly larger) compact subset in $\xi\in\R^n$.
\end{lem}

\Pf The assertion about the second condition is clear. To show that the first condition is also preserved, let points $x,y\in K$ be so close to each other that $||\sigma(x,\xi)-\sigma(y,\xi)||\le \delta$ for any $\xi$ (by the assumption of continuity of $\sigma$ in $x$ uniformly in $\xi$). We can write $\sigma(x,\psi_x(\xi))-\sigma(y,\psi_y(\xi))$ as a sum of two expressions: $(\sigma(x,\psi_x(\xi))-\sigma(x,\psi_y(\xi)))+(\sigma(x,\psi_y(\xi))-\sigma(y,\psi_y(\xi)))$. The norm of the second expression is estimated by $\delta$. 

The norm of the first expression is estimated as $\le ||d_\eta(\sigma(x,\eta))||\cdot ||\psi_x(\xi)-\psi_y(\xi)||$, where $\eta$ is some point on the segment joining $\psi_x(\xi)$ and $\psi_y(\xi)$, i.e. $\eta=t\psi_x(\xi)+(1-t)\psi_y(\xi)$ with $0\le t \le 1$. The first multiple is estimated as $\le C(1+||\eta||)^{-1}$. The second multiple is estimated as $\le ||\psi_x-\psi_y||\cdot ||\xi||$. Because $\psi_x$ is norm-continuous in $x$ and invertible for any $x\in K$, the map $\xi\mapsto \eta=t\psi_x(\xi)+(1-t)\psi_y(\xi)=t(\psi_x(\xi)-\psi_y(\xi))+\psi_y(\xi)$ is still an invertible linear map when $x$ and $y$ are close enough. Therefore $(1+||\eta||)^{-1}\le C_1(1+||\xi||)^{-1}$ for some constant $C_1$ depending on $x$ and $y$. This implies $||d_\eta(\hat\sigma(x,\eta))||\cdot ||\psi_x(\xi)-\psi_y(\xi)||\le C_2(1+||\xi||)^{-1}\cdot ||\psi_x-\psi_y||\cdot ||\xi||$, thus proving the uniform continuity of $\sigma(x,\psi_x(\xi))$. \cqfd

\medskip

Now we are ready to give a definition of symbols of `order $0$' and `negative order'. 
 
\begin{defi} Let $E$ be a (complex) vector bundle over $X$. A symbol $\sigma(x,\xi)$ of order $0$ is a bounded measurable section of the bundle $p^*(E)$ over $TX$ satisfying the following  conditions:

$1^{\rm o}$ For any compact subset in $x\in X$, $\sigma(x,\xi)$ is continuous in $x$ uniformly in $\xi$ outside of a compact subset in $\xi$. 

$2^{\rm o}$ For any compact subset in $x\in X$, $\sigma(x,\xi)$ is differentiable in $\xi$ outside of a compact subset in $\xi$, and for the exterior derivative $d_\xi$, there is an estimate for any compact subset $K\subset X$: $||d_\xi \sigma(x,\xi)||\le C\cdot (1+||\xi||)^{-1}$ with the constant $C$ which depends only on $K$ and $\sigma$ and not on $(x,\xi)$. 

We will say that a symbol $\sigma$ is of `negative order' if $||\sigma(x,\xi)||$ converges to $0$ uniformly in $x\in X$ on compact subsets of $X$ when $\xi\to\infty$. We will say that $\sigma$ is of `strongly negative order' if, additionally, $||\sigma(x,\xi)||$ converges to $0$ uniformly in $\xi$ when $x\to \infty$ in $X$. 
\end{defi}

\begin{rem} {\rm Condition $1^{\rm o}$ of definition 3.5 does not depend of the choice of local coordinates on the tangent bundle $T(X)$ in view of lemma 3.4. 

In fact, condition $2^{\rm o}$ of definition 3.5 was stated with this strong estimate on the first derivative of $\sigma$ because this estimate was required (for the change of coordinates) in the statement of lemma 3.4. However, there may be situations when, for example, there is a preferred choice for the trivialization of $T(X)$. In this case the precise estimate in condition $2^{\rm o}$ may be omitted. The minimal necessary requirement in condition $2^{\rm o}$ would be the vanishing of $d_\xi(\sigma)$ at infinity in $\xi$ uniformly on compact subsets in $x$.

Also note that the class of symbols of definition 3.5 includes the `classical' symbols (i.e. symbols homogeneous of order $0$ in the $\xi$ variable).} 
\end{rem} 

\subsection{Operators}

In our present approach, an operator corresponding to the symbol $\sigma$ will be constructed using the operator integration described in subsection 2.1. Let us fix a complex vector bundle $E$ over $X$ (equipped with a Hermitian metric). We will define operators acting on sections of $E$. 

Let $x\in X$ and $U_x\subset X$ is a small open ball of radius $r_x$ centered at $x$. Using the Riemannian exponential map $\exp_x$ (associated with the Levi-Civita connection), consider the open Euclidean ball $V_x$ of the same radius around the point $0\in T_x(X)$ as a Euclidean coordinate neighborhood for $U_x$. The tangent map $(\exp_x)_*:T(V_x)\to T(U_x)$ is almost an isometry when $r_x$ is sufficiently small (see \cite{Ka88}, proposition 4.3). In particular, $\exp_x$ defines an almost isometry $L^2(V_x)\to L^2(U_x)$. Also the family of first derivative linear maps $(\exp_x)_*$ is uniformly continuous in $\xi$ when $x$ varies and the point $p\in T_x(X)$ varies. This follows from the uniform dependence on the initial data for systems of linear differential equations (for Jacobi vector fields).

For any fixed $x\in X$, the symbol $\sigma_x(\xi)=\sigma(x,\xi):E_x\to E_x$ is a bounded function of $\xi\in T^*_x(X)$. If we denote the Fourier transform $L^2(T_x(X))\otimes E_x \to L^2(T^*_x(X))\otimes E_x$ by $\Phi$, then $ F_x=\Phi^{-1}\sigma_x\Phi$ is a bounded operator on $L^2(T_x(X))\otimes E_x$. This operator is just the operator with the symbol $\sigma_x(\xi)$. (Here $x$ is fixed, only $\xi$ varies.) 

Next, choose a function $\nu\in C_0^\infty([0,1))$ such that $0\le \nu \le 1$, $\nu(0)=1$, and $\nu(t)=0$ for $t\ge 1/2$. By theorem 3.1, $F_x$ commutes modulo compact operators with multiplication by functions from $C_c(T_x(X))$. In particular, it commutes modulo compact operators with the function $\nu_x(v)=\nu(||v||/r_x)$ (where we put $t=||v||/r_x$). 

Now we need a trivialization of $E$ over $U_x$. Of course, this trivialization will depend on $x$, so we need a continuous family of trivializations $\{E|_{U_x}\}$ over $X$. Using this and the almost isometry between $L^2(V_x)$ and $L^2(U_x)$ described above, we can transplant the operator $\nu_xF_x\nu_x$ to $L^2(E|_{U_x})$. Denote the resulting operator on $L^2(E|_{U_x})\subset L^2(E)$ by $F(x)$. 

Note that if the original symbol $\sigma$ is a symbol in the H\" ormander $\rho=1,\, \delta=0$ class, then each $F(x)$ belongs to the same class of PDOs in the H\" ormander theory and its symbol at the point $x$ is $\sigma(x,\xi)$ modulo symbols of lower order. This is essentially the homeomorphism invariance of the H\" ormander $\rho=1,\, \delta=0$ class (see \cite{Hor}, end of subsection 2.1).

Theorem 3.1 ensures that the operator $F(x)$ commutes modulo compact operators with the multiplication by functions from $C_0(X)$. (Each $F(x)$ commutes with $C_0(U_x)$, and because the cutting function $\nu_x$ has support inside $U_x$, the operator $F(x)$ commutes with all of $C_0(X)$.) 

Moreover, $F(x)$ is a norm-continuous function of $x\in X$ modulo compact operators. Indeed, by our assumptions, $\sigma(x,\xi)$ is continuous in $x$ uniformly in $\xi$. So the Euclidean operators $F_x(\xi)$ and $F_y(\xi)$ are norm-close to each other if $x$ and $y$ are close because they are just Fourier transforms of two uniformly close functions. Therefore the Euclidean operators $\nu_xF_x\nu_x$ and $\nu_yF_y\nu_y$ are also norm-close to each other. If $x$ and $y$ are close enough, both functions, $\nu_x$ and $\nu_y$, have their supports in $U_x\cap U_y$. The transition function between $E|_{U_x}$ and $E|_{U_y}$ will be close to the identity. Also the maps $L^2(V_x)\to L^2(U_x)$ and $L^2(V_y)\to L^2(U_y)$ will be norm-close to each other in $U_x\cap U_y$. Therefore, modulo compact operators, $F(x)$ and $F(y)$ will be norm-close to each other in $L^2(U_x\cap U_y)$.

Note that the concrete choice of the cutting functions $\nu_x$ has little  importance. Proposition 3.2, $2^{\rm o}$ shows that for two functions, $\nu_x$ and $\nu'_x$, the products $F(x)\nu_x$ and $F(x)\nu'_x$ are equal modulo compact operators if $\nu_x(0)=\nu'_x(0)$. 

For the operator integration procedure, we set $D=\K(L^2(E))$ and $\phi:C_0(X)\to \L(L^2(E))$ the multiplication by functions. The function $F$ constructed above will be considered as the function $X\to \M(D)/D$. The lifting of the operator integral $\int_X F(x)d\phi$ to $\L(L^2(E))$ will be denoted $\F$.

\begin{defi} We define `negative order' operators as those which become compact after multiplication by an element of $C_0(X)$. (This is the same as `locally compact' operators defined in 2.3.)
\end{defi}

In the case when there is a continuous, isometric, and proper action of a locally compact group $G$ on $X$ (and on the vector bundle $E$), and the symbol $\sigma$ is $G$-invariant, we can average $\F$ over $G$. For this, we take a cut-off function $\gtc$ on $X$, i.e. a non-negative function $\gtc\in C_b(X)$ whose support has compact intersection with any $G$-compact set in $X$ and which satisfies the property $\int_G g(\gtc)dg=1$. Let $\sum_i \beta_i=1$ be a partition of unity on $X/G$ so that all $\beta_i\ge 0$, all supports $\supp\beta_i$ are compact in $X/G$ and form a locally finite covering of $X/G$. We will consider all functions $\beta_i$ as $G$-invariant functions on $X$. 

Using remark 2.2 and proposition 2.8, we take $Av_G (\beta_i\gtc\F)$ for every $i$. The operators $\beta_i\gtc\F$ have compact support because all $\beta_i\gtc$ have compact support and $\F$ has proper support - see remark 2.2 and theorem 2.4, $5^{\rm o}$. After that we redefine $\F$ as $\sum_i Av_G (\beta_i\gtc\F)$. Since the original operator $\F$ had support in an $\epsilon$-neighborhood of the diagonal of $X\times X$ for some $\epsilon>0$ (see the proof of 2.4, $5^{\rm o}$), the averaged operator $\F$ will have the same property. Therefore it will be properly supported. It will differ from the original $\F$ by an operator of `negative order'.

\begin{theo} The correspondence between the symbol $\sigma$ and the operator $\F$ constructed out of it has the following properties (modulo negative order symbols and negative order operators):

$1^{\rm o}$ Composition of symbols $\mapsto$ composition of operators.

$2^{\rm o}$ $\sigma^* \mapsto \F^*$.

$3^{\rm o}$ If the symbol $\sigma$ is bounded at infinity in $\xi$ by $C>0$, i.e. for any compact subset $K\subset X$ and any $x\in K$, $\limsup_{\xi\to \infty} ||\sigma(x,\xi)||\le C$, then the norm of the operator $\F$ in $\L(L^2(E))/\K(L^2(E))$ does not exceed $C$.

$4^{\rm o}$ The operator corresponding to a negative order symbol has negative order, and the operator corresponding to a strongly negative order symbol is compact.
\end{theo}

\Pf All assertions follow directly from the construction and theorem 2.4. To obtain the last assertion concerning strongly negative order symbols, one can either use corollary 2.6, or one may approximate the symbol by symbols with compact support in $x\in X$. \cqfd

\begin{theo} Let $\sigma(x,D)$ be a bounded, properly supported pseudo-differen\-tial operator of order $0$ with  the symbol $\sigma$ in the H\" ormander $\rho=1,\,\delta=0$ calculus, and let $\F$ be the operator constructed out of $\sigma$ above. Then, the difference $\sigma(x,D)-\F$ is an operator of `negative order'. 

In particular, the operator $\F$ has a `principal symbol', i.e. a symbol modulo symbols of negative order.
\end{theo}

\Pf Both $\sigma(x,D)$ and $\F$ are bounded. We need to prove that $\sigma(x,D)$ and $\F$ coincide modulo $\K$ on any compact piece of $X$. Let $\{U_i\}$ be a locally finite covering of $X$ consisting of Riemannian balls $U_{x_i}$ and $\sum_i \alpha_i^2=1$ the corresponding partition of unity. Let us present $\sigma(x,D)$ (modulo $\K$) as the sum: $\sum_i \alpha_i\sigma(x,D)\alpha_i$. If diameters of the balls $U_{x_i}$ are small enough (on our compact piece of $X$), this sum is equal to $\sum_i \alpha_i\nu_{x_i}\sigma(x,D)\nu_{x_i}\alpha_i$ (because $\nu_{x_i}$ is equal to $1$ in a certain neighborhood of the point $x_i$).

Let us compare the latter sum with the integral sum $\sum_i \alpha_i F(x_i)\alpha_i$. For any $i$, both $\alpha_i F(x_i)\alpha_i$ and $\alpha_i\sigma(x,D)\alpha_i$ are PDOs of order $0$ with compactly supported distributional kernel. The difference between their symbols goes to $0$ uniformly in $i$ when the radii of the balls $U_i$ go to $0$ (on our compact piece of $X$). Now the assertion follows from the usual norm estimate results: lemma 2.5 and \cite{Hor}, corollary 2.2.3. 

Concerning the definition of a principal symbol for the operators $\F$, we define it as the principal symbol of $\sigma(x,D)$. The `principal symbol' depends only on the operator modulo operators of negative order (see \cite{Hor}, end of subsection 2.1). 
\cqfd

\medskip

Finally, a remark about notation. In the previous construction we have actually used not the symbol of an operator but its Fourier transform.

\begin{defi} We will denote the tangent bundle $T(X)$ by $\tau$ and the tangent space $T_x(X)$ at $x\in X$ by $\tau_x$. The algebra $C^*(\tau_x)$ is the \Cst-algebra of the abelian group $\tau_x$ with its Euclidean topology. The algebra of continuous sections (vanishing at infinity) of the field $\{C^*(\tau_x),x\in X\}$ will be denoted $C^*_\tau(X)$. The scalar symbols $\sigma$ of order $0$ live in $C_b(TX)$, their Fourier transforms (cosymbols) live in $\M(C^*_\tau(X))$.
\end{defi}

\medskip

We end up this section with a simple technical lemma which will be needed later in the proof of our index theorems. 

Let $Y\to X$ be a smooth, locally trivial fiber bundle, both $X$ and $Y$ are complete Riemannian manifolds. We will denote the fiber over a point $x\in X$ by $Y_x$. Assume that the transition functions are isometries in the Riemannian metrics on all $Y_x$. Let $\E$ be the Hilbert module over $C_0(X)$ defined by the continuous field $L^2(Y_x),\, x\in X$. Consider the Hilbert space $H=\E\otimes_{C_0(X)} L^2(X)$. For any $e\in \E$, define the operator $\theta_e\in \L(L^2(X),H)$ by $\theta_e(l)=e\otimes l\in H$ for any $l\in L^2(X)$. 

Now a small reminder: Suppose we have a bounded operator $F\in \L(L^2(X))$. Then a \emph{$K$-theoretic connection} (cf. \cite{CS}, A1; \cite{Ka88}, 2.6) for the operator $F$ is an operator $\tilde F\in \L(H)$ such that for any $e\in \E$, the operators $\theta_e F-\tilde F \theta_e$ and $\theta_e F^*-\tilde F^* \theta_e$ are compact. (In fact, when there are some additional graded vector bundles over $X$ and $Y$ and we consider graded $L^2$-spaces, then the multiple $(-1)^{deg(e)\cdot deg(F)}$ has to be introduced in front of $\tilde F$ and $\tilde F^*$, so that the commutators become graded).

\begin{lem} In the notation above, suppose that $F$ is a PDO with the symbol $\sigma$. Let us not distinguish between tangent and cotangent bundles (using the Riemannian metrics) and lift $\sigma$ from $TX$ to $TY$ via the projection $Y\to X$. Denote the lifted symbol by $\tilde \sigma$. Then the PDO $\tilde F$ on the space $H$ with the symbol $\tilde \sigma$ (constructed by means of our coarse PDO calculus) is a $K$-theoretic connection for $F$.
\end{lem}

\Pf If $Y$ is a direct product $X\times Y_x$, the assertion is obviously true. But the lifted symbol commutes with the transition functions. \cqfd

\section{\Cst-algebras associated with a group action}

Let $X$ be a locally compact, $\sigma$-compact space equipped with a continuous proper action of a locally compact group $G$. The left-invariant Haar measure on $G$ will be denoted $dg$, and the modular function of $G$ by $\mu$. The stability subgroup of any point $x\in X$ will be denoted $M_x$. All stability subgroups are compact because the action is proper. We will assume that the Haar measure of any compact group is normalized to the total mass $1$. 

For any Hilbert $C_0(X)$-module $E$, the fiber of $E$ over $x\in X$ will be denoted $E_x$. The space of compactly supported elements of $E$, i.e. $E\cdot C_c(X)$, will be denoted $E_c$ (or $C_c(E)$ if we consider $E$ as a continuous field of vector spaces).

For a proper $G$-algebra $B$, we will use the notation $B^G$ for a certain subalgebra of $\M(B)$ consisting of $G$-invariant elements - see \cite{Ka88}, definition 3.2.
  
We start with the construction which originates in the work of M. Rieffel. Here we will follow \cite{KS03}, section 5. 

\begin{defi} Let $B$ be a $G-C_0(X)$-algebra, where $X$ is a proper $G$-space, and $E$ a Hilbert $G-B$-module. Denoting by $C_c(X)$ the set of compactly supported functions on $X$, let $E_c=C_c(X)\cdot E$. We define a pre-Hilbert module structure over $C_c(G,B)$ on $E_c$ by 
$$e\cdot b=\int_G g(e)\cdot g(b(g^{-1}))\cdot \mu(g)^{-1/2}dg\in E_c,$$
$$(e_1,e_2)(g)=\mu(g)^{-1/2}(e_1,g(e_2))_E \in C_c(G,B),$$
where $e,e_1,e_2\in E_c,\; b\in C_c(G,B)$. We will denote by $\E$ the Hilbert $C^*(G,B)$-module which is the completion of $E_c$ in the norm defined by the above $C_c(G,B)$ inner product on $E_c$.
\end{defi}

\begin{rem} {\rm The positivity of the inner product comes from the embedding $i:\E\subset C^*(G,E)$ (onto a direct summand) given by $i(e)(g)=\mu(g)^{-1/2}c^{1/2}g(e)$, where $c$ is a cut-off function on $X$. The Hilbert module $C^*(G,E)$ is defined in \cite{Ka88}, 3.8. It is isomorphic to $E\otimes_B C^*(G,B)$.}  
\end{rem}

\begin{lem} $\K(\E)\simeq \K(E)^G$. In particular, when $E=B$, we have $\K(\E)=B^G$, and when $E=B=C_0(X)$, we get $\K(\E)=C_0(X/G)$. 
\end{lem}

\Pf The `rank 1' operators $\theta_{e_1,e_2}$ in $\K(\E)$ (i.e. $\theta_{e_1,e_2}(e)=e_1(e_2,e)$) are easily seen to be $\theta_{e_1,e_2}(e)=\int_G \theta_{g(e_1),g(e_2)}(e)dg,$ 
where $e_1,e_2,e\in E_c$, and on the right side, the notation $\theta$ is used for the `rank 1' operators in $\K(E)$. This proves the assertion. \cqfd

\begin{lem} When $E=L^2(G,B)$, then $\E\simeq C^*(G,B)$ as a Hilbert module over $C^*(G,B)$. 
\end{lem}

\Pf Consider $L^2(G,B)$ with the right $G$-action: $g(e)(t)=\mu(g)^{1/2}g(e(tg))$ instead of the usual left one. (These two $G$-actions correspond to each other under the automorphism $e(g)\mapsto \mu(g)^{-1/2} e(g^{-1})$ of $L^2(G,B)$.) Let us use elements $e=e(g)\in E_c$ with compact support in $g\in G$ (i.e. elements from $C_c(X)\cdot C_c(G,B)$). Then it is easy to check that the isomorphism $\E\simeq C^*(G,B)$ is given by the formula: $e(g)\mapsto g(e(g))$. \cqfd

\medskip

The last thing that we need to mention before we state the result that we need is the canonical unitary representation $u:G\to \M(C^*(G,B))$ given on elements $b\in C_c(G,B)$ by $u_h(b)(g)=h(b(h^{-1}g))$ for any $g,h\in G$. 

Now lemmas 4.3 and 4.4 imply the following result (hidden in the proof of theorem 5.4 of \cite{KS03}), which is a generalization of the well known similar fact for compact groups (cf. e.g. \cite{GHT}, 11.2):

\begin{theo} $C^*(G,B)\simeq \K(L^2(G,B))^G$, where we use the right $G$-action on $L^2(G,B)$. Under this isomorphism, the canonical representation $u:G\to \M(C^*(G,B))$ transforms into the left translation action of $G$ on $L^2(G,B)$. 
\end{theo}

\Pf The isomorphism follows directly from 4.3 and 4.4. For the second statement, take $E=L^2(G,B)$. Then the left translation action of $G$ on $E$ exactly corresponds to the representation $u$ under the map $e(g)\mapsto g(e(g))$ given in the proof of lemma 4.4: $u_h(ge(g))=h(h^{-1}g(e(h^{-1}g))=g(e(h^{-1}g))$, which is equal to $g(e_1(g))$ for $e_1(g)=e(h^{-1}g)$. \cqfd

\medskip

Next, we will prove a version of the structure theorem of P. Green (\cite{Gr80}, 2.13). Let us assume that $X$ is a fibered product $X=G\times_M Z$, where $M$ is a compact subgroup of $G$ and $Z$ an $M$-space. There is a natural projection $\pi:X\to G/M$ in this case, and $\pi(Z)=(M)\in G/M$. Let $B_Z$ be the part of $B$ over $Z$, i.e. $B_Z=B/(B\cdot C_0(X-Z))$. (In the Riemannian manifold case which will be considered later, $Z$ will be a small Euclidean disk orthogonal to the orbit $O_x$ at the point $x\in X$, and $M$ will be the stability subgroup at $x$.) 

In this case, Green's result is that $C^*(G,B)\simeq C^*(M,B_Z)\otimes \K(L^2(G/M))$. We need only a weak version of it, and we want to avoid using a measurable cross-section in the construction.

It is clear that $B$ is the algebra defined by the continuous field of \Cst-algebras $\{B_t, t\in G/M\}$, where $B_t=t(B_Z)$. This means that $B$ is isomorphic to the algebra of continuous sections (vanishing at infinity) of the fiber bundle of algebras $G\times_M B_Z \to G/M$ associated with the principal bundle $G\to G/M$. We will use the notation $(C_0(G)\otimes B_Z)^M$ for this algebra.

\begin{prop}  In the above assumptions and notation, if $X=G\times_M Z$, then
$C^*(G,B)\simeq (B_Z\otimes \K(L^2(G)))^M$, where $M$ acts on $B_Z\otimes\K(L^2(G))$ diagonally (with the action on $L^2(G)$ by right translations). Under this isomorphism, the canonical representation $u:G\to \M(C^*(G,B))$ transforms into the left translation action of $G$ on $L^2(G)$.
\end{prop} 

\Pf By theorem 4.5, $C^*(G,B)\simeq (B\otimes \K(L^2(G)))^G$. This can be rewritten as $((C_0(G)\otimes B_Z)^M\otimes \K(L^2(G)))^G$. Here $G$ acts on $C_0(G)$ by left translations and $M$ by right translations, and these two actions commute. The action of $G$ on $L^2(G)$ is by right translations, and the action of $M$ on $L^2(G)$ is trivial. Therefore $C^*(G,B)$ is isomorphic to $(C_0(G)\otimes B_Z\otimes \K(L^2(G)))^{M\times G}$. 

The isomorphism stated in the proposition comes from the obvious isomorphism: $(C_0(G)\otimes A)^G\simeq A$ for $A=B_Z\otimes \K(L^2(G))$. Note that the inverse of the latter isomorphism maps $a\in A$ to $g\mapsto \{g(a)\}$ in $C_b(G,A)$. Because $M$ acts on $C_b(G)$ by right translations, we get the action of $M$ on $C_b(G,A)$: $\{g(a)\}\mapsto \{gm(a)\}$, which corresponds to the diagonal action $a\mapsto m(a)$ of $M\subset M\times G$ on $A$. 

The last assertion follows from the similar assertion in theorem 4.5. \cqfd

\section{Leaf-wise PDOs}

In this section, we give a definition of symbol algebras for leaf-wise PDOs, followed by the construction of leaf-wise operators. We also consider a more geometric construction for some classical geometric operators like a leaf-wise Dirac operator.  

The notation and assumptions concerning the manifold $X$ will be the same as stated in section 3: $X$ will be a complete Riemannian manifold equipped with a proper isometric action of a Lie group $G$, $T(X)$ and $T^*(X)$ the tangent and cotangent bundles of $X$ respectively (isomorphic via the Riemannian metric). The tangent manifold will be denoted $TX$, and the projection $TX\to X$ by $p$.

\subsection{Notation and preliminaries.}

Consider the map $f:G\times X\to X, f(g,x)=(g(x),x)$ defining the $G$-action on $X$. For each $x\in X$, let $f_x:G\to X$ be the map $f$ restricted to $G\times \{x\}$, i.e. $f_x(g)=g(x)$. The orbit $\O_x$ through $x$ is the image of $f_x$. We have: $L^2(G/M_x)\subset L^2(G)$. The image of $L^2(G/M_x)$ consists of all $L^2$-functions on $G$ right-invariant under $M_x$. However, although $\O_x\simeq G/M_x$, the $L^2$-spaces of $\O_x$ and $G/M_x$ are not in general isometrically isomorphic because their metric (and measure) may be quite different. We will use two fiber bundles: the trivial one $q:G\times X\to X$ with the projection $q$ onto the second multiple and also $p:T(X)\to X$.

Recall some notation from \cite{Ka16}, sections 6 and 7. We denote by $\gtg$ the Lie algebra of the group $G$. Let $f'_x: \gtg\to T_x(X)$ be the tangent map (first derivative of $f_x$) at the identity of $G$, and $f'^*_x:T_x^*(X)\to \gtg^*$ the dual map. It is easy to see that for any $x\in X,\; g\in G,\; v\in \gtg$, one has: $g(f'_x(v))=f'_{g(x)}(Ad(g)(v))$. 

Let us consider the trivial vector bundle $\gtg_X=X\times \gtg$ over $X$ with the $G$-action given by $(x,v)\mapsto (g(x),Ad(g)(v))$. Because the $G$-action on $X$ is proper, there exists a $G$-invariant Riemannian metric on $\gtg_X$. Equivalently, one can say that there exists a smooth map from $X$ to the space of Euclidean norms on $\gtg$: $x\mapsto ||\cdot ||_x$,  such that for any $x\in X,\; g\in G,\; v\in \gtg$, one has: $||Ad(g)(v)||_{g(x)}=||v||_x$. 

We denote by $f':\gtg_X\to T(X)$, the map defined by $(x,v)\mapsto f'_x(v)$ at any $x\in X$. From the above formulas, it is clear that this map is $G$-equivariant. Note that by multiplying our Riemannian metric on $\gtg_X$ by a certain strictly positive $G$-invariant function, we can also arrange that the following condition is satisfied: for any $v\in \gtg,\; ||f'_x(v)||\le ||v||_x$. We will assume that this is the case, and so we have $||f'_x||\le 1$ for any $x\in X$.

Let us identify $\gtg_X$ with its dual bundle via the Riemannian metric. The fiber of $\gtg_X$ at the point $x$ will be denoted $\gtg_x$. Consider a continuous field of subspaces of the tangent bundle $T(X)$ defined as $\{\im f'_x\subset T_x(X),\;x\in X\}$. Call this field $\Gamma$. This is the field of tangent spaces to the orbits of $G$. Another field is defined as $\{\im f'^*_x\subset \gtg^*_x,\; x\in X\}$. These two fields are isometrically isomorphic. The isomorphism is given by the isometry $(f'f'^*)^{-1/2}f':\im f'^*\to \im f'$. (\emph{Note that we do not assume that $f'$ is invertible.}) The continuous sections of these two fields are $\{f'(v),\; v\in C_0(\gtg_X)\}$ and $\{f'^*(v),\; v\in C_0(T(X))\}$. \emph{We will denote by $\Gamma$ any of these two fields of vector spaces and by $\Gamma_x$ the fiber of $\Gamma$ at $x\in X$.} 

We define a $G$-invariant map $\{\varphi_x=(f'_xf'^*_x)^{1/2},\;x\in X:T^*(X)\to T^*(X)\}$. (\emph{Note the change in notation from \cite{Ka16}, section 6, where the same notation $\varphi_x$ was used for the map $f'_xf'^*_x$.}) We also define a $G$-invariant quadratic form $q$ on covectors $\xi\in T_x^*(X)$: $q_x(\xi)=(f'_xf'^*_x(\xi),\xi)=||f'^*_x(\xi)||_x^2=||\varphi_x(\xi)||^2$. Our previous assumption $||f'_x||\le 1$ implies  that $q_x(\xi)\le ||\xi||^2$ for any covector $\xi$. Being restricted to any orbit, the form $q_x(\xi)$ can serve as a Riemannian metric on this orbit. This is compatible with the natural topology on each orbit and on the orbit space globally. 

Note that a covector $\xi\in T_x^*(X)$ is orthogonal to the orbit passing through $x$ if and only if $q_x(\xi)=0$. The subspace $\{(x,\xi), q_x(\xi)=0\}\subset TX$ is $T_GX$ in Atiyah's notation. 

The Lie algebra of the stability subgroup $M_x$ will be denoted $\gtm_x$, and its orthogonal complement by $\gtm_x^\perp\subset \gtg_x^*$. Note that $\gtg_x/\gtm_x$ is isomorphic (via $f'_x$) to the tangent space $T_x(O_x)$, and $\gtm_x^\perp=(\gtg_x/\gtm_x)^*$ is isomorphic (via $f'^*_x$) to the cotangent space $T^*_x(O_x)$.

\subsection{Symbol algebras}

For any $x\in X$, both vector spaces $\gtg_x$ and $\tau_x=T_x(X)$ can be considered as abelian groups, so there are the corresponding convolutional \Cst-algebras: $C^*(\gtg_x)$ and $C^*(\tau_x)$. They are isomorphic via the Fourier transform to $C_0(\gtg_x^*)$ and $C_0(\tau_x^*)$ respectively. The map $(f'^*_x)^*:C_b(\gtg_x^*)\to C_b(\tau_x^*)$ corresponds to the natural map $(f'_x)_*:\M(C^*(\gtg_x))\to \M(C^*(\tau_x))$ via the Fourier transform.

Let us denote the algebras of continuous sections (vanishing at infinity of $X$) of the fields $\{C^*(\gtg_x),\, x\in X\}$ and $\{C^*(\tau_x),\, x\in X\}$ by $C^*_\gtg(X)$ and $C^*_\tau(X)$ respectively. We will denote by $\gtG X$ the total space of the fiber bundle $\gtg_X\simeq X\times \gtg$ and by $q$ its projection to $X$. Then $\{C_0(\gtg_x^*),\, x\in X\}$ are exactly the fibers of $C_0(\gtG X)$ over $X$. Similarly, the fibers of the algebra $C_0(TX)$ over $X$ form a continuous field of algebras $\{C_0(\tau^*_x),\, x\in X\}$. Therefore $C_0(\gtG X)\simeq C^*_\gtg(X)$ and $C_0(TX)\simeq C^*_\tau(X)$ via the Fourier transform. 

The important conclusion is that there are natural maps $(f'^*)^*:C_b(\gtG X)\to C_b(TX)$ and $(f')_*:\M(C^*_\gtg(X))\to \M(C^*_\tau(X))$ related via the Fourier transform. 

\medskip

Now we can define the symbol algebra for leaf-wise operators. This will be the algebra of operators of `negative order'.  

\begin{defi} Using the maps $C_b(\gtG X)\to C_b(TX)$ and $\M(C^*_\gtg(X))\to \M(C^*_\tau(X))$, we define the symbol algebra $\Slf(X)$ for the `negative order' leaf-wise operators as the image of $C_0(\gtG X)$ in $C_b(TX)$. This is a $C_0(X)$-algebra. The convolutional symbol algebra $C^*_\Gamma(X)$ (isomorphic to $\gtS_{lf}(X)$ via the Fourier transform) is defined as the image of $C^*_\gtg(X)$ in $\M(C^*_\tau(X))$.

For the Hilbert module $E$ over $C_0(X)$ of continuous sections of a (finite-dimensional) vector bundle over $X$, we consider $\gtS_{lf}(E)=\K(E\otimes_{C_0(X)} \Slf(X))$ as the corresponding symbol algebra associated with $E$. 
\end{defi}

At this point, it will be useful to make some additional clarification concerning the algebra $\Slf(X)$. Note that the map $f'_x:\gtg_x\to\tau_x$ factors as follows: $\gtg_x\to\gtg_x/\gtm_x\simeq T_x(O_x)\subset\tau_x$. The dual map $f'^*_x$ factors as $\tau^*_x\to T^*_x(O_x)\simeq (\gtg_x/\gtm_x)^*\simeq \gtm_x^\perp \subset \gtg_x^*$. 

\begin{defi} We will consider $\{C_0(\gtm_x^\perp),\,x\in X\}$ as a continuous field of algebras. The continuous sections of this field are, by definition, the restrictions of continuous sections belonging to $C_0(\gtG X)$.
We define $C_0((\gtG/\gtM) X)$ as the algebra of continuous sections (vanishing at infinity of $X$) of the continuous field $\{C_0(\gtm_x^\perp),\,x\in X\}$. 
\end{defi}

Because the map $C_0(T^*_x(O_x))\to C_b(\tau^*_x)$ induced by the restriction $\tau^*_x\to T^*_x(O_x)$ is injective, we get the following simple fact:

\begin{prop} The algebra $\Slf(X)$ is isomorphic to $C_0((\gtG/\gtM) X)$. \cqfd
\end{prop}

\begin{ass} Before discussing actual leaf-wise symbols and operators, here are the two options in which these objects will be defined:

$1^{\rm o}$ The $G$-equivariant option: symbols and operators are $G$-invariant. 

$2^{\rm o}$ The stabilizer-invariant option: For any open subset $U_x=G\times_{M_x}Z_x\subset X$, the symbol restricted to the slice $Z_x$ is invariant with respect to the action of the stability subgroup $M_x$.
\end{ass}

Since the map $C_0(\gtG X)\to \gtS_{lf}(X)$ is surjective, it extends to the  surjective map of the multiplier algebras: $\M(C_0(\gtG X))\to \M(\gtS_{lf}(X))$ (\cite{APT}, theorem 4.2).

We will use the product structure of the bundle $\gtG X=X\times \gtg$. 

\begin{defi} An element $\sigma(x,\xi)\in \L(\gtS_{lf}(E))$ (where $\xi\in \tau_x$) will be called a leaf-wise symbol of order $0$ if it is an image of an element $\gtb(x,\eta)\in \L(E\otimes_{C_0(X)} C_0(\gtG X))$ (where $\eta\in \gtg_x$) satisfying the following assumptions:

$1^{\rm o}$ $\gtb(x,\eta)$ is continuous in $x$ uniformly in $\eta$ (on compact subsets of $X$); 

$2^{\rm o}$ $\gtb(x,\eta)$ is differentiable in $\eta$, and its exterior derivative $d_\eta(\gtb)$ vanishes at infinity in $\eta$ uniformly in $\eta$ on compact subsets of $X$. 

A symbol $\sigma$ is of `negative order' if it becomes an element of $\K(\gtS_{lf}(E))$ after being multiplied by any element of $C_0(X)$.
\end{defi}

\begin{rem} {\rm The strong estimate on $\gtb(x,\eta)$ required in definition 3.5, $2^{\rm o}$, is not needed in the item $2^{\rm o}$ of definition 5.5 because of an almost canonical choice of coordinate systems in the tangent spaces along the orbit (and the stabilizer-invariance of the symbols).}
\end{rem}

\begin{rem} {\rm Definition 5.5 was given with the assumption that all our symbols are continuous functions. Recall from definition 3.5 that one can also use symbols which are not everywhere continuous, e.g. `classical' symbols. In order to take this possibility into account, definition 5.5 should be modified. This refers both to $\sigma(x,\xi)$ and $\gtb(x,\eta)$. The terminology that has to be used is given below (cf. also definition 3.5).}
\end{rem}

\begin{term}We will say that a function $a(x,\eta)$ satisfies a certain condition outside of a compact set in $\eta$ if for any compact set $L'$ in the variable $x$ there is a compact set $L''$ in the variable $\eta$ such that the said condition holds for all $(x,\eta)$ with $x\in L',\eta\notin L''$.
\end{term}

\subsection{Clifford algebras and leaf-wise Dirac type operators}

Our goal in this and the next subsection is to construct a multiplier of the algebra $C^*(G,\K(E))$ corresponding to a leaf-wise operator. The Hilbert module $E$ here is the module of continuous sections (vanishing at infinity) of a (complex) vector bundle over $X$. In the present subsection we use a specific way to do it for some geometric operators of Dirac type. First recall some notation from \cite{Ka16}. 

\medskip

Let $\V$ be a real vector bundle over $X$ equipped with a $G$-invariant Riemannian metric and ${\it Cliff}(\V,Q)$ the Clifford algebra bundle associated with the quadratic form $Q(v)=||v||^2$ on $\V$. We denote by $Cl_\V(X)$ the complexification of the algebra of continuous sections of ${\it Cliff}(\V,Q)$ over $X$, vanishing at infinity of $X$. With the sup-norm on sections, this is a \Cst-algebra. It will be denoted $Cl_\V(X)$. When $\V=\tau=T^*(X)$, this algebra will be denoted $\Ct(X)$, and when $\V=\gtg_X$, this algebra will be denoted $Cl_\gtg(X)$. Recall that according to \cite{Ka16}, theorem 2.7, the algebras $\Ct(X)$ and $C_0(TX)$ are $KK$-equivalent.

Let $f':\gtg_X\to \tau=T(X)$ be the tangent map defined in subsection 5.1 and $f'^*:\tau^*\to \gtg_X^*$ the dual map. The \Cst-subalgebra of $\Ct(X)$ generated by scalar functions $C_0(X)\subset \Ct(X)$ together with the subspace $\{f'(v)\}$ for all continuous sections $v$ of $\gtg_X$ vanishing at infinity will be denoted  $Cl_\Gamma(X)$. A similar \Cst-subalgebra of $Cl_\gtg(X)$ generated by scalar functions $C_0(X)\subset Cl_\gtg(X)$ together with the subspace $\{f'^*(w)\}$ for all covector fields $w\in \tau^*$ vanishing at infinity is isomorphic to $Cl_\Gamma(X)$ because the spaces of continuous sections $\{f'(v),\; v\in C_0(\gtg_X)\}$ and $\{f'^*(w),\; w\in C_0(\tau^*)\}$ are isometrically isomorphic (see subsection 5.1). (We will not distinguish between these two \Cst-algebras).

Note that on each orbit $\O$ in $X$, the algebra $\CG(X)$ restricted to $\O$ is isomorphic to $\Ct(\O)$ where the Riemannian metric on $\O$ is defined by the family of quadratic form $q_x, x\in \O$, because $q_x(\xi)=||f'^*_x(\xi)||_x^2$ (see subsection 5.1).

\medskip

\noindent{\bf \large Construction.}

\medskip

The usual classical Dirac type operators are the Dirac operator acting on a spin bundle, the Euler characteristic operator, the signature operator, and the Dolbeault operator (in the almost complex manifold case). The construction is similar in all these cases. Here we will give the details for the construction of the Euler characteristic operator. We will call it `Dirac operator' because it will serve as the basis for the construction of the important Dirac element later in 8.8. The notation $\Lambda^*(X)$ in this construction will mean the \emph{complexified} vector bundle $\Lambda^*(X)$. We will also denote by $\Lambda^*(X)$ the Hilbert module of its continuous sections vanishing at infinity. (All PDOs that we consider act on \emph{complex} vector bundles.)

\medskip

For each orbit $\O_x\simeq G/M_x$, we will first construct a canonical differential Dirac operator acting on sections of the vector bundle $\Lambda^*(\O_x)$. Since the metric on the orbit (induced from $X$) is different from the metric on $G/M_x$, we start with an operator on $G/M_x$. The space of $L^2$-sections of the vector bundle $\Lambda^*(G/M_x)$ is isomorphic to $(L^2(G)\otimes \Lambda^*(\gtm_x^\perp))^{M_x}$, where $M_x$ acts on $L^2(G)$ by right translations. 

For any $v\in \gtg$, let $g_v(t)$ be the one-parameter subgroup of $G$ corresponding to the vector $v$. The infinitesimal right translation by $g_v(t)$ on $L^2(G)$ defines the left-invariant differential operator $\partial/\partial v$ on smooth elements of $L^2(G)$. One can also apply it on $C_c^\infty(G)\otimes \Lambda^*(\gtm_x^\perp)$ using the composition of the (infinitesimal) right translation on $L^2(G)$ and the usual left $\gtg$-action on $\Lambda^*(\gtm_x^\perp)$. 
Note that when $v\in \gtm_x$, this differential operator is $0$ on $C_c^\infty(G/M_x)$ and on $(C_c^\infty(G)\otimes \Lambda^*(\gtm_x^\perp))^{M_x}$. (We consider $C^\infty(G/M_x)$ as a subspace of $C^\infty(G)$, $L^2(G/M_x)$ as a subspace of $L^2(G)$, and $\Lambda^*(\gtm_x^\perp)$ as a subspace of $\Lambda^*(\gtg^*)$.)

Let us choose any basis $\{\tilde v_k\}$ in $\gtm_x^\perp$ and the dual basis $\{v_k\}$ in $\gtg/\gtm_x$. Define the differential Dirac operator on $(C_c^\infty(G)\otimes \Lambda^*(\gtm_x^\perp))^{M_x}$ by the formula: 
$$\D_{G/M_x}=-i\sum_k (\ext(\tilde v_k)+\inter(\tilde v_k))\partial/\partial v_k.$$
This definition does not depend on the choice of the basis and gives an unbounded $G$-invariant operator on $(C_c^\infty(G)\otimes \Lambda^*(\gtm_x^\perp))^{M_x}$, and an unbounded, essentially self-adjoint operator on $(L^2(G)\otimes \Lambda^*(\gtm_x^\perp))^{M_x}$. Indeed, one can check that this operator preserves the $M_x$-invariants in $C_c^\infty(G)\otimes \Lambda^*(\gtm_x^\perp)$. The $G$-invariance is also clear because the $G$-action is defined by left translations on $L^2(G)$. 

\begin{rem}
{\rm The important property of the operator $\D_{G/M_x}$ on $G/M_x$ is that it is a natural restriction of the corresponding operator $\D_G$ on $G$ which acts on $C_c^\infty(G)\otimes \Lambda^*(\gtg^*)$. More generally, if $M_y\subset M_x$, then $\D_{G/M_x}$ is the restriction of $\D_{G/M_y}$. (We repeat that  we consider $C^\infty(G/M_x)$ as a subspace of $C^\infty(G)$, $L^2(G/M_x)$ as a subspace of $L^2(G)$, and $\Lambda^*(\gtm_x^\perp)$ as a subspace of $\Lambda^*(\gtg^*)$.)}
\end{rem}

All operators $\D_{G/M_x}$ are elliptic on the corresponding spaces $G/M_x$ and their principal symbols are easy to calculate: $\sigma_{\D_{G/M_x}}=\ext(\eta)+\inter(\eta)$. The principal symbol of $\D_{G/M_x}^2$ is $||\eta||^2$, i.e. the same as the principal symbol of the Laplace operator $\Delta_{G/M_x}$. Let us put $\F_{G/M_x}=\D_{G/M_x}(\D_{G/M_x}^2+1)^{-1/2}$. This is a bounded self-adjoint operator on $L^2(\Lambda^*(G/M_x))$, and $1-\F_{G/M_x}^2=(\D_{G/M_x}^2+1)^{-1/2}$ is an operator which becomes compact after being multiplied by any continuous function vanishing at infinity.

Now we can descend (push forward) all these operators to the orbits in $X$ using the homeomorphisms $f_x:G/M_x\simeq \O_x$. The operators $\D$ and $\F$ on the orbit $\O$ will be denoted $\D_\O$ and $\F_\O$ respectively. They depend only on the orbit, not on the point $x$ of the orbit. 

Using the embedding $\gtm_x^\perp\simeq T_x(\O_x)\subset T_x(X)=\tau_x$, one can redefine the operators $\D_{\O_x}$ as operators acting on the bundles $\Lambda^*(X)|_{\O_x}$. The formula for the new operator $\D_{\O_x}$ will be the same as before. This will give the operators $\D_\O$ and $\F_\O$ on sections of the bundle $\Lambda^*(X)$ over the orbit $\O$. 

The constructed family of operators $\{\F_\O\}$ is a multiplier of the algebra $C^*(G,\K(\Lambda^*(X)))$. In fact, for any tubular neighborhood $U=G\times_{M_x}Z_x$ of an orbit $\O_x\subset X$, we have: $C^*(G,\K(\Lambda^*(U))\simeq \K(L^2(G)\otimes \Lambda^*(X)|_{Z_x})^{M_x}$. By remark 5.9, all operators of our family are restrictions of one operator from $G$ to all orbits, so the (strict) continuity of the family $\{\F_\O\}$ is clear. Moreover, the family $\{f\cdot(1-\F_\O^2)\}$ is norm-continuous for any $f\in C_0(X)$ and gives an element of $C^*(G,\K(\Lambda^*(X)))$.

The principal symbol of the descended operator $\F_{\O_x}$ (calculated using the change of variables formula in the H\" ormander calculus) is equal $\sigma_\F(x,\xi)=(\ext(f'^*_x(\xi))+\inter(f'^*_x(\xi)))(1+||f'^*_x(\xi)||^2)^{-1/2}=(\ext(f'^*_x(\xi))+\inter(f'^*_x(\xi)))(1+q_x(\xi))^{-1/2}$. The operator $1-\F_{\O_x}^2$ has principal symbol $(1+q_x(\xi))^{-1}$.

Note that $f'^*_x(\xi)$ is actually an element of $\gtg_x^*$. If we want our symbol to be expressed in terms of covectors on $X$, we need to apply the isomorphism of the subsection 5.1 between the implementations of the field $\Gamma$ on $\gtG X$ and on $TX$ provided by the isometry $(f'f'^*)^{-1/2}f':\im f'^*\to \im f'$. Applying this isometry, we get $(f'_xf'^*_x)^{-1/2}(f'_xf'^*_x)(\xi))=\varphi_x(\xi)$ instead of $f'^*_x(\xi)$ in the above formulas. Therefore, $\sigma_{\D_{O_x}}(x,\xi)=(\ext(\varphi_x(\xi))+\inter(\varphi_x(\xi)))$, and the principal symbol of the operator $\F_{O_x}$ is $(\ext(\varphi_x(\xi))+\inter(\varphi_x(\xi)))(1+q_x(\xi))^{-1/2}$.

\begin{notat} {\rm It is convenient to shorten the notation $\ext(\varphi_x(\xi))+\inter(\varphi_x(\xi))$ to $c(\varphi_x(\xi))$, where $c$ means the usual Clifford action on the exterior algebra. We will use the notation $c$ for Clifford multiplication in the rest of the paper.}
\end{notat}

\begin{defi} The family $\{\F_\O\}$ will be called the leaf-wise Dirac operator on $\Lambda^*(X)$ and denoted $\gtD_\Gamma$.
\end{defi}

We have proved the following 

\begin{prop} There exists a leaf-wise Dirac operator $\gtD_\Gamma$ on $\Lambda^*(X)$. The principal symbol of 
$\gtD_\Gamma$ is $c(\varphi_x(\xi))(1+q_x(\xi))^{-1/2}$. This Dirac operator defines an element of the algebra $\M(C^*(G,\K(\Lambda^*(X))))$. Furthermore, $C_0(X)\cdot (1-\gtD_\Gamma^2)\subset C^*(G,\K(\Lambda^*(X)))$. \cqfd
\end{prop}

\begin{rem} {\rm The constructed operator $\gtD_\Gamma$ can be considered as an element of the group $KK^G(C_0(X), C^*(G,\K(\Lambda^*(X)))$, but this is only one of our goals here. In fact, based on the operator $\gtD_\Gamma$, we will construct an element of a different $KK$-group later, in example 8.8. That will be the basic leaf-wise Dirac element that we need for our main results.}
\end{rem}

\subsection{General construction of leaf-wise operators}

This subsection contains a construction of an element of $\M(C^*(G,\K(E)))$ corresponding to a leaf-wise operator with the symbol satisfying the assumptions of definition 5.5.

Given a leaf-wise symbol $\sigma$, we choose the corresponding element $\gtb(x,\eta)\in \L(E\otimes_{C_0(X)} C_0(\gtG X))$ (definition 5.5). Actually we will regard $\sigma$ as an element of $\L(E\otimes_{C_0(X)} C^*_\Gamma(X))$ (cf. definition 5.1) and $\gtb$ as an element  of $\L(E\otimes_{C_0(X)} C^*_\gtg(X))$. The condition $1^{\rm o}$ imposed on $\gtb$ in definition 5.5 means that $\gtb$ is norm-continuous in $x\in X$.

We will use the fact that the space $X$ can be covered with the fibered products of the form: $U_x=G\times_{M_x} Z_x$, where $M_x$ is the stability subgroup at the point $x\in X$ and $Z_x$ is a small open Euclidean disk. According to proposition 4.6, the \Cst-algebra $C^*(G, \K(E|_{U_x}))$ is isomorphic to $\K(L^2(G)\otimes E|_{Z_x})^{M_x}$, where $M_x$ acts on $L^2(G)\otimes E|_{Z_x}$ diagonally (with the action on $L^2(G)$ by right translations).

For each fibered product $U_x=G\times_{M_x} Z_x$, one gets the corresponding family of pseudo-differential operators on $L^2(G)$ parametrized by points of $Z_x$ (according to the construction of section 3). Because the family of symbols is assumed to be stabilizer-invariant (by assumption 5.4), this family of PDOs naturally descends to a multiplier of $\K(L^2(G)\otimes E|_{Z_x})^{M_x}$ by averaging over $M_x$ (see subsection 2.2). Using a partition of unity on $X/G$, one gets a multiplier of $C^*(G,\K(E))$. In the case of $G$-invariant symbols, the averaging over $G$ gives a $G$-invariant PDO.

This construction does not depend on a particular choice of the element $\gtb$ used in it. The operators on $L^2(G)$ constructed using $\gtb$ do depend on $\gtb$, but when they are pushed forward to the orbits, this dependence vanishes. The pushforward map is the map $f$ defined at the beginning of subsection 5.1. It transforms the $\gtb$-symbol into the $\sigma$-symbol.   

\begin{rem} {\rm The construction of the Dirac operator (subsection 5.3) can also be performed using the method described here.}
\end{rem}

Finally, let us define `compact operators' for $\M(C^*(G,\K(E)))$ as elements of $C^*(G,\K(E))$. \emph{We reiterate definition 3.7 of `negative order' operators as those which become compact after being multiplied by a function from $C_0(X)$.} 

\begin{lem} If the symbol is of `negative order', then the corresponding operator is of `negative order'.
\end{lem}

\Pf The algebra $C_0(\gtG X)$ has the function $(1+||\eta||^2)^{-1}$ as its fiber-wise (for each $x\in X$) strictly positive element. The image of this element in $\gtS_{lf}(X)$ is the element $(1+q_x(\xi))^{-1}$ for each $x\in X$. In subsection 5.3, we have already proved our assertion for the operator $\{(1+\D_{O_x}^2)^{-1}\}$ with such symbol. It remains to note that the operator integration method preserves the relation of order, modulo compact operators (see theorem 3.8). Therefore, if our assertion is true for the operator with the symbol $(1+q_x(\xi))^{-1}$, it is true for the operator with any `negative order' symbol. \cqfd

\section{Poincare duality}
 
 We will need a certain generalization of the Poincare duality \cite{Ka16}, theorem 4.6. The proof of this result is exactly the same as the proof of \cite{Ka16}, 4.6. A sketch of the proof will be given below.
 
 \smallskip
 
Let us start with the local Bott element. Let $U$ be a small tubular neighborhood of the diagonal $\Delta\subset X\times X$ (homeomorphic to $TX$) such that for any $(x,y)\in U$ there is a minimal geodesic connecting $x$ and $y$. The local Bott element $[\B_X]\in \cR KK^G(X;C_0(X), C_0(U)\otimes_{C_0(X)} C_0(TX))$ (defined in \cite{Ka16}, 2.9) is a family of point-wise Bott elements $\{\beta_x\}, x\in X$. Each $\beta_x$ is the Bott element in $K_0(C_0(U_x)\otimes C_0(T_x(X)))$, where $U_x=U\cap (x\times X)$. Embedding $U$ into $X\times X$ we get $[\B_X]\in \cR KK^G(X;C_0(X),C_0(X)\otimes C_0(TX))$.

If $A$ is a $C_0(X)$-algebra, we can tensor the local Bott element with $1_A$ over $C_0(X)$. We get an element $[\B_A]\in \cR KK(X;A, A\otimes C_0(TX))$.  

\begin{theo} Let $X$ be a complete Riemannian manifold with an isometric proper action of a locally compact group $G$. For any separable $G-C_0(X)$-algebra $A$ and a separable $G$-algebra $B$,
$$KK^G_*(A,B)\simeq {\cal R}KK^G_*(X;A,C_0(TX)\otimes B)$$
$$\simeq \cR KK^G_*(X;C_0(TX)\otimes_{C_0(X)}A,C_0(X)\otimes B).$$ 
\end{theo}

\noindent\emph{Sketch of proof.}  To deal with the first isomorphism, we define the homomorphisms:
$$\mu=[\B_A]\otimes_{A}:KK_*^G(A,B)\to \cR KK_*^G(X;A,C_0(TX)\otimes B)$$
and $$\nu=\otimes_{C_0(TX)} [\D_X]: \cR KK_*^G(X;A,C_0(TX)\otimes B)\to KK_*^G(A,B),$$
where $[\D_X]\in K^0_G(C_0(TX))$ is the Dolbeault element (\cite{Ka16}, 2.8). The associativity of the product and theorem 2.10 \cite {Ka16} easily imply that $\nu\cdot \mu=id$. In the opposite direction, the proof is exactly the same as in \cite{Ka16}, 4.6.

The second isomorphism comes from the Bott periodicity. The groups $\cR KK^G_*(X;C_0(TX)\otimes_{C_0(X)}A,C_0(X)\otimes B)$ and $\cR KK^G_*(X;A,C_0(TX)\otimes B)$ are isomorphic: take an $\cR KK$-product (over $C_0(X)$) with $1_{C_0(TX)}$ and use the $\cR KK$-equivalence (Bott periodicity) $C_0(TX)\otimes_{C_0(X)}C_0(TX)\simeq_{KK} C_0(X)$. \cqfd

\begin{rem} {\rm It is useful to have Poincare duality also for manifolds with boundary. We will indicate the necessary changes. Let $\bar X$ be a manifold with boundary $\partial X$ and interior $X$. To construct $[\B_{\bar X}]\in \cR KK^G(\bar X;C_0(\bar X), C_0(\bar X)\otimes C_0(TX))$, we can extend $X$ to a complete Riemannian manifold $Y=\bar X\cup_{\partial X} (\partial X\times [0,\infty))$ and then restrict the element $[\B_Y]$ to $\bar X$ using the fact that $X$ and $Y$ are diffeomorphic. 

Now $A$ will be a $C_0(\bar X)$-algebra. We tensor $[\B_{\bar X}]$ with $1_A$ over $C_0(\bar X)$ and get an element $[\B_A]\in \cR KK^G(\bar X;A, A\otimes C_0(TX))$. The Dolbeault element does not change.
The resulting Poincare duality isomorphism will be the following: $KK^G_*(A,B)\simeq {\cal R}KK^G_*(\bar X;A,C_0(TX)\otimes B)$.}
\end{rem}

\section{K-theory of symbol algebras}

Now we come to the symbol algebras. The symbol algebra $\gtS_{lf}(X)$ for leaf-wise `negative order' operators was defined in 5.1 as the image of $C_0(\gtG X)$ in $C_b(TX)$ under the map $f'^*:TX\to \gtG X$. For any $x\in X$, the strictly positive element of $C_0(\gtG X)$ in the fiber over $x$ is $(1+||\eta||^2)^{-1}$. Therefore the algebra $\gtS_{lf}(X)$ is fiberwise generated (inside $C_b(TX)$) by the field of elements $(1+||f'^*_x(\xi)||^2)^{-1}=(1+q_x(\xi))^{-1}$, with an additional condition that elements of $\gtS_{lf}(X)$ vanish at infininity in $x\in X$ uniformly in $\xi$. 

Recall now a similar algebra for transversal PDOs. It was defined in \cite{Ka16} (definitions 6.2, 6.3) and denoted $\gtS_G(X)$ and $\gtS_\Gamma(X)$. Here we change slightly this definition according to our more flexible approach to PDOs. We also prefer to change slightly the notation for this algebra. 

\begin{defi} The symbol algebra $\gtS_{tr}(X)$ for `negative order' transversal PDOs is a subalgebra of the algebra of bounded functions $b\in C_b(TX)$ satisfying the conditions:

$1^{\rm o}$ $b(x,\xi)$ is continuous in $x$ uniformly in $\xi$ and vanishes (uniformly in $\xi$) when $x\to\infty$. 

$2^{\rm o}$ $b(x,\xi)$ is differentiable in $\xi$, and for the exterior derivative in $\xi$, there is an estimate: $||d_\xi(b)||\le C\cdot (1+||\xi||)^{-1}$ for any compact subset $K\subset X$, with the constant $C$ which depends only on $b$ and $K$ and not on $(x,\xi)$. 

$3^{\rm o}$ For any $\varepsilon>0$, there exists $c>0$ such that for any $x\in X,\;\xi\in T^*_x(X)$,
$|b(x,\xi)|\le c\cdot (1+q_x(\xi))(1+||\xi||^2)^{-1}+\varepsilon.$

Similarly defined, symbols of operators on a finite-dimensional vector bundle $E$ over $X$ form a Hilbert module $\gtS_{tr}(E)$ over $\gtS_{tr}(X)$.
\end{defi}

\begin{rem} {\rm It is clear from condition $3^{\rm o}$ that the field of elements $(1+q_x(\xi))(1+||\xi||^2)^{-1}$ is the field of strictly positive elements in the fibers of the algebra $\gtS_{tr}(X)$ over $X$. Because the products $(1+q_x(\xi))^{-1}\cdot (1+q_x(\xi))(1+||\xi||^2)^{-1}=(1+||\xi||^2)^{-1}$ form a field of strictly positive elements in the fibers of the algebra $C_0(TX)$ over $X$, we have: $\gtS_{lf}(X)\cdot \gtS_{tr}(X)= C_0(TX)$. In particular, there is a natural map of $C_0(X)$-algebras: $\gtS_{lf}(X)\otimes_{C_0(X)} \gtS_{tr}(X)\to C_0(TX)$.}
\end{rem}

\begin{defi-lem} Denote the Hilbert $C_0(X)$-module corresponding to the field of Hilbert spaces $\{L^2(\Lambda^*(\tau_x))\}$ over $X$ by $\T$. There is the Bott-Dirac operator on each space $L^2(\Lambda^*(\tau_x))$, namely, $\cC_x=D_x(1+D_x^2)^{-1/2}$, where $D_x=(d_\xi+d^*_\xi+\ext(\xi)+\inter(\xi))$. Denote the corresponding operator on $\T$ by $\cC$. The homomorphism $\psi:\CG(X)\otimes_{C_0(X)}\gtS_{lf}(X)\to \L(\T)$ is just multiplication on the $\gtS_{lf}(X)$ part, and on the Clifford part of $\CG(X)$ it is given on (real) covectors by $\xi\mapsto \ext(\xi)+\inter(\xi)$. The triple $(\T,\psi, \cC)$ defines a Dirac type element of the group $\cR KK^G(X;\CG(X)\otimes_{C_0(X)}\gtS_{lf}(X),C_0(X))$ which will be denoted $[d_{lf}]$.
\end{defi-lem}

\Pf We need only to prove the compactness of commutators between $\cC$ and the image of $\psi$. Any element $a\in \psi(\CG(X)\otimes_{C_0(X)}\gtS_{lf}(X))$ can be considered as a function of $x\in X$ and $\xi\in \tau_x^*$. We can assume that $a$ is a smooth covector field. For any $x\in X$, the direction of the covector $a_x(\xi)$ belongs to the subspace $\Gamma_x$ tangent to the orbit passing through $x$. Moreover, because $a$ belongs to the image of $\psi$, as a function of $\xi$ (with $x$ frozen), $a_x(\xi)$ is constant in $\xi$ in the directions orthogonal to $\Gamma_x$, and we can assume that in the direction of $\Gamma_x$ it has compact support in $\xi$ (since we can assume that an element of $\Slf(X)$ is an image of a compactly supported element of $C_0(\gtG X)$).

It is clear from the above that the graded commutator $[a_x(\xi),\ext(\xi)+\inter(\xi)]$, which is equal to the scalar product $2(a_x(\xi),\xi)$, vanishes for any $\xi$ orthogonal to $\Gamma_x$. On the other hand, for $\xi$ in the direction of $\Gamma_x$, the above commutator is bounded because we assumed that $a_x(\xi)$ has compact support in that direction. Therefore, the commutator $[a_x(\xi),\ext(\xi)+\inter(\xi)]$ is always bounded uniformly in $x$. 

Since $a_x(\xi)$ has compact support in the direction of $\Gamma_x$ and is constant in the direction orthogonal to $\Gamma_x$, all first $\xi$-derivatives of $a_x$ vanish at infinity in $\xi$. Therefore the (graded) commutator $[a_x,d_\xi+d_\xi^*]$ is bounded on $L^2(\Lambda^*(\tau_x))$. 

Now we will use the method of lemma 4.2 \cite{Ka16}. Write $\cC_x=(2/\pi)\int_0^\infty D_x(1+\lambda^2+D_x^2)^{-1}d\lambda$. For any $a_x$ as above, the (graded) commutator $[a_x,D_x]$ is bounded on $L^2(\Lambda^*(\tau_x))$ and depends norm-continuously on $x\in X$. The operators $(1+\lambda^2+D_x^2)^{-1}$ and $D_x(1+\lambda^2+D_x^2)^{-1}$ are compact on $L^2(\Lambda^*(\tau_x))$ and depend norm-continuously on $x\in X$. The commutator $[a_x, \cC_x]$ is presented as the integral 
$$(2/\pi)\int_0^\infty (1+\lambda^2+D_x^2)^{-1}((1+\lambda^2)[a_x,D_x]+D_x[a_x,D_x]D_x)(1+\lambda^2+D_x^2)^{-1}d\lambda$$
which converges uniformly in $\lambda$, so the result is compact and depends norm-continuously on $x\in X$. \cqfd 

\begin{defi} Two separable $G-C_0(X)$-algebras $D_1$ and $D_2$ will be called $\cR KK$-dual if there are elements 
$$\alpha\in \cR KK^G(X;D_1\hotimes_{C_0(X)} D_2,C_0(X)), \;\; \beta\in \cR KK^G(X;C_0(X),D_1\hotimes_{C_0(X)} D_2)$$ 
such that the product with these two elements gives isomorphisms:
$$\cR KK^G(X; A\hotimes_{C_0(X)} D_1,B) \simeq \cR KK^G(X; A,B\hotimes_{C_0(X)} D_2)$$
and
$$\cR KK^G(X; A\hotimes_{C_0(X)} D_2,B) \simeq \cR KK^G(X; A,B\hotimes_{C_0(X)} D_1)$$
for any two separable $G-C_0(X)$-algebras $A$ and $B$. 
\end{defi}

\begin{rem} {\rm Assume that the above elements $\beta$ and $\alpha$ have their product over ${D_1\hotimes_{C_0(X)}D_2}$ equal to $1_{C_0(X)}\in \cR KK^G(X; C_0(X),C_0(X))$. Then the $\cR KK$-duality of $D_1$ and $D_2$ can be proved using the same kind of rotation homotopy that is used in Atiyah's proof of the Bott periodicity. Both algebras $D_1$ and $D_2$ must allow this rotation homotopy. 

Indeed, for example, taking the product of $a\in \cR KK^G(X; A,B\hotimes_{C_0(X)} D_2)$ first with $1_{D_1}$ and then with $1_{D_2}$, we arrive at an element of the group
$$\cR KK^G(X; A\hotimes_{C_0(X)} D_1\hotimes_{C_0(X)} D_2,B\hotimes_{C_0(X)} D_2\hotimes_{C_0(X)} D_1\hotimes_{C_0(X)} D_2).$$ 
Then we need to take the product with $\beta$ and $\alpha$. We would like the result to be equal to the initial element $a$. This would have been the case if we took the product with $\alpha$ over the last couple of $D_1\hotimes_{C_0(X)} D_2$. But in fact it is the previous couple $D_2\hotimes_{C_0(X)} D_1$ that is required to be used for the product. So we need a rotation homotopy which interchanges the two copies of $D_2$ in the last argument of the $\cR KK$-group.}
\end{rem}

\begin{theo} The algebras $\gtS_{lf}(X)$ and $\CG(X)$ are $\cR KK$-dual. The algebra $\CG(X)\otimes_{C_0(X)}\gtS_{lf}(X)$ is $\cR KK$-equivalent to $C_0(X)$.
\end{theo}

\Pf The element $[d_{lf}]$ of definition 7.3 will serve as $\alpha$. The element of  $\cR KK^G(X;C_0(X),\CG(X)\otimes_{C_0(X)}\gtS_{lf}(X))$ which will be denoted $[\gtf_\Gamma]$, and will serve as $\beta$, is defined as follows: 

The Hilbert module is $\CG(X)\otimes_{C_0(X)}\gtS_{lf}(X)$. The operator $\gtf_\Gamma$ is defined by the covector field $\{\gtf_{\Gamma_x}(\xi)=\varphi_x(\xi)(1+q_x(\xi))^{-1/2},\;x\in X\}$.

Let $D=\CG(X)\otimes_{C_0(X)}\gtS_{lf}(X)$. To show that the product $\beta\otimes_D \alpha$ is equal to $1_{C_0(X)}$, we need to produce a homotopy of the product operator $N_1\gtf_\Gamma+N_2\cC$ to $\cC$. (The operators $N_1,N_2$ are those used in the $KK$-product construction.)

First, let us look at the graded commutator $Q=[N_1\gtf_\Gamma+N_2\cC,\cC]$. Note that $[\varphi_x(\xi),d_\xi+d^*_\xi]$ is constant in $\xi$ and $[\varphi_x(\xi),\ext(\xi)+\inter(\xi)]=2(\varphi_x(\xi),\xi)\ge 0$. Using the same integral presentation as in lemma 7.3, both for $\gtf_\Gamma$ and $\cC$, we will use the same calculation of the commutator $[\gtf_\Gamma,\cC$] as in the proof of theorem 8.9 \cite{Ka16}. This calculation and the reasoning given in the proof of 8.9 \cite{Ka16} show that $[\gtf_\Gamma,\cC]$ is non-negative modulo compact operators, so the same is true for $Q$. Now the existence of the homotopy joining $N_1\gtf_\Gamma+N_2\cC$ to $\cC$ follows from lemma 11 \cite{Sk84}.

In the remaining part of the proof we will use a rotation homotopy. For each $x\in X$, the fiber of $\gtS_{lf}(X)$ over $x$ is the image of $C_0(\gtg_x)$ in $C_b(\tau_x)$.  So we will use the rotation of $\gtg_x\times \gtg_x$. Similarly for $\CG(X)$, where the fiber over $x$ is just the corresponding Clifford algebra.

Let us first show that $\alpha\otimes_{C_0(X)}\beta=1_D$. This is very similar to the proof of the Bott periodicity in \cite{Ka80}, theorem 5.7. This product can be written as $(\phi:D\to \L(\T\otimes_{C_0(X)} D), \cC\;\kkprod\; \gtf_\Gamma)$. The map $\phi$ factors through $\M(D\hotimes D)$. There is an obvious rotation homotopy of $D\hotimes D$, as already explained. The operator $\gtf_\Gamma$ rotates together with $\phi$. The operator $\cC$ does not rotate. At the end of this rotation homotopy we arrive at the product $(\beta\otimes_D \alpha)\otimes_{C_0(X)} 1_D=1_D\in \cR KK^G(D,D)$. The details on choosing the homotopy of operators $N_1,N_2$ involved in the $KK$-product construction see \cite{Ka80}, the proof of theorem 5.7.

To complete the proof, it remains to show the existence of a rotation homotopy indicated in the remark 7.5. Let $D_1=\CG(X), D_2=\gtS_{lf}(X)$, and consider the rotation for the $\cR KK$-group in the displayed formula in the remark 7.5. If $\E$ is the Hilbert module for $a\in \cR KK^G(X; A,B\hotimes_{C_0(X)} D_2)$, we need to compose the homomorphism 
$$A\hotimes_{C_0(X)} D_1\hotimes_{C_0(X)} D_2 \to \L(\E\hotimes_{B\hotimes_{C_0(X)}D_2} B\hotimes_{C_0(X)}D_2\hotimes_{C_0(X)}D_1\hotimes_{C_0(X)} D_2)$$
with the rotation interchanging the last two copies of $D_2$. Clearly, there is a map of $D_2\hotimes_{C_0(X)} D_2$ into the right hand side, and we can compose it with the rotation of the map $D_2\to \M(D_2\hotimes_{C_0(X)} D_2)$, as already explained. \cqfd

\medskip

Before we move further, we need to correct an error in \cite{Ka16} at the end of definition 7.1, where it was claimed that the algebra $\CG(X)\hotimes_{C_0(X)}\CG(X)$ is isomorphic to the algebra of compact operators on the Hilbert module $C_0(\Lambda^*_\Gamma(X))$ over $C_0(X)$. Unfortunately, this is a wrong statement, which was used a number of times in sections 7 and 8 of \cite{Ka16}. The algebras $\CG(X)\hotimes_{C_0(X)}\CG(X)$ and $C_0(X)$ are not Morita equivalent in general, although they are Morita equivalent in a special case when $\Gamma$ is defined as the field of continuous sections of a vector bundle. But for a general $\Gamma$ we will use the $\cR KK$-duality between $\CG(X)$ and $\Slf(X)$  provided by the previous theorem:

\begin{cor} If $A$ and $B$ are any separable $G-C_0(X)$-algebras, then 
$$\cR KK^G(X; A\hotimes_{C_0(X)}\CG(X),B)\simeq \cR KK^G(X; A,B\hotimes_{C_0(X)}\Slf(X)),$$
$$\cR KK^G(X; A\hotimes_{C_0(X)}\Slf(X),B)\simeq\cR KK^G(X; A,B\hotimes_{C_0(X)}\CG(X)).$$
 \end{cor}

Combined with the Poincare duality 6.1, this gives a corrected version of theorem 7.8 \cite{Ka16} -- see theorem 7.11 below. Here is another corollary of Poincare duality 6.1 and the $\cR KK$-duality:

\begin{prop} Let $G$ be a compact group and $X$ a compact $G$-manifold. If the $G-C(X)$-algebras $A$ and $A^{\rm dual}$ are $\cR KK$-dual then 
$$K^*_G(A)\simeq K_*^G(A^{\rm dual}\hotimes_{C_0(X)} C_0(TX)).$$
\end{prop}

\Pf The left side is isomorphic to $\cR KK_*^G(X;A, C_0(TX))$ by theorem 6.1, which in its turn is isomorphic to $\cR KK_*^G(X;C_0(X),A^{\rm dual}\hotimes_{C_0(X)} C_0(TX))$ by $\cR KK$-duality, and coincides with the right side in the statement because $X$ is compact. \cqfd
 
 \begin{notat} {\rm We simplify the notation $Cl_\Gamma(X)\otimes_{C_0(X)} C_0(TX)$ to $\CG(TX)$. (Note that the latter \emph{is not} related with the orbits of $G$ on the manifold $TX$.) We will also denote the algebra $\Ct(X)\hotimes_{C_0(X)}\CG(X)$ by $\CtG(X)$.}
\end{notat}

Now we are going to restate and reprove theorem 7.4 of \cite{Ka16} which relates the symbol algebra $\gtS_{tr}(X)$ with the algebra $\CG(TX)$. We remind that the symbol algebra $\gtS_{tr}(X)$ is not separable, so to state the $KK$-equivalence one needs a separable version of the $KK$-theory introduced by G. Skandalis in \cite{Sk88}, pp. 571-572 (cf. also the beginning of section 7 in \cite{Ka16}). 

\medskip

Here are some definitions from \cite{Ka16} that we will need: 

\medskip

First, a family of Dirac operators on tangent spaces of $X$ (\cite{Ka16}, 2.5). For any $x\in X$, let $H_x=L^2(T_x^*(X))\otimes Cl_{\tau_x}$ be the Hilbert module  over $Cl_{\tau_x}$. The Dirac operator on this Hilbert module is defined by $\gtd_x=\sum_{k=1}^{\dim X} (-i)c(e_k)\partial/\partial \xi_k$ in any orthonormal basis $\{e_k\}$ of $T^*_x(X)$, where $c(e_k)$ are the left Clifford multiplication operators and $\xi_k$ are the coordinates of $T^*_x(X)$ in the basis $\{e_k\}$. The family of operators $\{\gtd_x/(1+\gtd_x^2)^{-1/2}\}$ on the family of Hilbert modules $\{H_x,\,x\in X\}$ defines an element $\cR KK^G(X; C_0(TX),\Ct(X))$ denoted $[d_\xi]$. 

Next, the fiber-wise Bott element  
$[\B_{\xi,\Gamma}]\in \cR KK^G(X;\CtG(X),\Str(X))$ (\cite{Ka16}, 7.2), which is given by the pair $(\Str(E),\beta)$, where $E=\Lambda^*(X)$ is considered as a Hilbert module over $C_0(X)$, $\Str(E)=E\otimes_{C_0(X)}\Str(X)$, the left action of $\CtG(X)$ on $\Str(E)$ is defined on (real) covectors by 
$$\xi_1\oplus\xi_2\mapsto i(\ext(\xi_1)-\inter(\xi_1))+(\ext(\xi_2)+\inter(\xi_2)).$$
Here $\xi_1$ is a section of $\tau=T(X)$ and $\xi_2$ a section of $\Gamma \subset \tau$. The operator $\beta$ is the Bott operator defined by $\beta(\xi)=(\ext (\xi)+\inter (\xi))/(1+||\xi||^2)^{1/2}$.

If we forget about $\Gamma$ in this definition of $[\B_{\xi,\Gamma}]$ and replace $\Str(X)$ with $C_0(TX)$ we get the tangent Bott element $[\B_\xi]\in \cR KK^G(X;\Ct(X),C_0(TX))$ (\cite{Ka16}, 2.6).

The operator $\gtf_\Gamma=\{\varphi_x(\xi)(1+q_x(\xi))^{-1/2},\;x\in X\}$ will be considered here as an element of $\M(Cl_\Gamma(TX))$. The element $[\gtf_\Gamma]\in \cR KK^G(X;\gtS_{tr}(X),Cl_\Gamma(TX))$ is defined by the natural homomorphism $\gtS_{tr}(X)\to \M(Cl_\Gamma(TX))$ together with the operator $\gtf_\Gamma$.

\begin{theo} The algebras $\CG(TX)$ and $\CtG(X)$ are $\cR KK^G$-equivalent, and $\gtS_{tr}(X)$ is $\cR KK^G_{sep}$-equivalent to each of them.
\end{theo}

\Pf The $\cR KK^G$-equivalence between the algebras $\CG(TX)$ and $\CtG(X)$ follows from the $\cR KK^G$-equivalence between $C_0(TX)$ and $\Ct(X)$ (theorem 2.7 \cite{Ka16}) by tensoring it over $C_0(X)$ (i.e. $\otimes_{C_0(X)}$) with $\CG(X)$. This $\cR KK^G$-equivalence is given in fact by the product with the elements $[d_\xi]$ and $[\B_\xi]$ defined above which are inverses of each other (\cite{Ka16}, theorem 2.7). 

The product $[\B_{\xi,\Gamma}]\otimes_{\Str(X)}[\gtf_\Gamma]$ is given by the pair 
$$(\psi:\CtG(X)\to \L(E\hotimes_{C_0(X)}\CG(TX)), N_1(\beta\hotimes 1)+N_2(1\hotimes \gtf_\Gamma)),$$
where $E=\Lambda^*(X)$. Here $N_1,N_2$ are the operators entering the $KK$-product construction: $N_1\cdot \Str(X)\subset C_0(TX)$ and $N_2\cdot \Slf(X)\subset C_0(TX)$. 

Let us make a rotation homotopy of this element. Consider the homomorphism $\CG(X)\hotimes_{C_0(X)}\CG(X)\to \L(E\otimes_{C_0(X)}\CG(TX))$, where the first copy of $\CG(X)$ maps according to the map $\psi$ and the second copy maps into $\M(\CG(TX))$. We will rotate this homomorphism by using the rotation of $\Gamma_x\times\Gamma_x$ ($x\in X$). The result of this homotopy will be $\gtf_\Gamma$ acting on $E$ and $\CG(X)$ acting on $\CG(TX)$. The operator will turn into $N_1(\beta\hotimes 1)+N_2(\gtf_\Gamma\hotimes 1)$. 

Finally, observe that this operator is homotopic to $\beta\hotimes 1$. This follows from lemma 11 \cite{Sk84} after we use the same argument as in the proof of theorem 7.6 above based on the fact that the graded commutator $[N_1(\beta\hotimes 1)+N_2(\gtf_\Gamma\hotimes 1),\beta\hotimes 1]\ge 0$ modulo compact operators. (In fact, $[\gtf_\Gamma(\xi),\beta(\xi)]$ differs by a positive multiple from the scalar product $(\varphi_x(\xi),\xi)$.) 

As a result, we get the equality: $[\B_{\xi,\Gamma}]\otimes_{\Str(X)}[\gtf_\Gamma]=[\B_\xi]\otimes_{C_0(X)} 1_{\CG(X)}$. Because $[\B_\xi]$ is an invertible element, we have proved that both $[\B_{\xi,\Gamma}]$ and $[\gtf_\Gamma]$ are one-sided invertible. To show that these two elements are actually invertible, we will prove that the product $[\gtf_\Gamma]\otimes_{C_0(X)}([d_\xi]\otimes_{C_0(X)}1_{\CG(X)})\otimes_{C_0(X)}[\B_{\xi,\Gamma}]$ is equal $1_{\Str(X)}$.

Let us first look at the product $([d_\xi]\otimes_{C_0(X)}1_{\CG(X)})\otimes_{C_0(X)}[\B_{\xi,\Gamma}]$. We will denote it $[d_{tr}]$. It is given by the pair $(\pi:\CG(TX)\to \L(\T\hotimes_{C_0(X)} \Str(X)), F)$, where the Hilbert $C_0(X)$-module $\T$ is defined in lemma 7.3. The map $\pi$ of $\CG(TX)$ is given by multiplication by functions of $C_0(\tau_x)$ on $L^2(\tau_x)$ (for each tangent space $\tau_x,\,x\in X$), and on (real) covectors of $\CG(TX)$ the map $\pi$ is $\xi\mapsto \ext(\xi)+\inter(\xi)$. The operator $F=N_1(A\hotimes 1)+N_2(1\hotimes \beta)$, where $A=(d_\xi+d_\xi^*)(1+(d_\xi+d_\xi^*)^2)^{-1/2}$, both $d_\xi+d_\xi^*$ and $\beta$ use the same exterior algebra variables, but $d_\xi+d_\xi^*$ acts on the $L^2(\tau_x)$ spaces, whereas $\beta(\xi)$ is the multiplication operator on $\Str(X)$.

The product $[\gtf_\Gamma]\otimes_{\CG(TX)} [d_{tr}]$ is given by the pair $(\tilde\pi:\gtS_{tr}(X)\to \L(\T\otimes_{C_0(X)}\gtS_{tr}(X)), N_1(A\hotimes 1)+N_2(1\hotimes (\beta+\gtf_\Gamma))$. Here $\Str(X)$ acts on $\T$ and $\gtf_\Gamma$ is the multiplication operator on $\Str(X)$. 

We first deform the latter product using the operator homotopy $N_1(A\hotimes 1)+N_2(1\hotimes (\beta+t\gtf_\Gamma)),\; 1\ge t\ge 0$. After that we make a rotation homotopy using the fact that $\tilde\pi$ factors through the homomorphism $\M(\Str(X)\otimes_{C_0(X)} \Str(X)) \to \L(\T\otimes_{C_0(X)}\gtS_{tr}(X))$. As a result of this rotation homotopy, the map $\tilde\pi$ transforms into the identity map $\Str(X)\to\Str(X)$, and $\beta$ transforms into the multiplication operator which acts on the $L^2(\tau_x)$  spaces. Therefore the operator will be homotopic to the operator $\cC$ of lemma 7.3. This means that our product turns into $1_{\Str(X)}$. \cqfd
 
 \medskip

We collect the most important results of this section in the next theorem:

\begin{theo}
 $$\cR KK^G_*(X;C_0(X),\Slf(X))\simeq \cR KK^G_*(X;\CG(X),C_0(X))$$
 $$\simeq \cR KK^G_*(X;\CtG(X),\Ct(X))\simeq K^*_G(\CtG(X)).$$
 
$$\cR KK^G_*(X;C_0(X),\gtS_{tr}(X))\simeq \cR KK^G_*(X;C_0(X),\CtG(X))$$
$$\simeq \cR KK^G_*(X;\Slf(X),\Ct(X))
\simeq K^*_G(\Slf(X)).$$
\end{theo}

\Pf The first isomorphism comes from corollary 7.7. The next isomorphism is the Clifford periodicity provided by the fact that $\tau$ is a vector bundle. The final isomorphism in the second line comes from the Poincare duality theorem 6.1. In the third line, the first isomorphism is given by theorem 7.10. The next isomorphism comes from corollary 7.7. The last isomorphism of the fourth line is provided by the Poincare duality 6.1.
\cqfd

\begin{cor} When $G$ and $X$ are compact we get:

\;\;\; $K^G_*(\Str(X))\simeq K_G^*(\Slf(X))$ and $K^G_*(\Slf(X))\simeq K_{G,sep}^*(\Str(X))$.
\end{cor}

\section{Elliptic symbols and index groups}

We start with a definition of certain special $KK$-groups related with crossed products. Then we give an overview of symbol and index groups for various types of operators. In the last part of this section, there are some important results on interrelations between different index and symbol groups.

\subsection{Special $KK$-groups for crossed products}

Let $B$ be a $G$-algebra. Consider the canonical unitary representation $u:G\to\M(C^*(G,B)): g\mapsto u_g$ (see subsection 4.1 before theorem 4.5). Then the $G$-action on $C^*(G,B)$ is given by $g(b)=u_gbu_g^{-1}$, where $b\in C^*(G,B)$. For any Hilbert module $\E$ over $C^*(G,B)$, one can define a structure of a $G-C^*(G,B)$-module as follows.  For $e\in \E$, put $g(e)=e\cdot u_g^{-1}$. An easy check shows that this definition of the $G$-action on $\E$ satisfies the axioms of a $G-C^*(G,B)$-module (\cite{KaJOT}, definition 2.1, $5^{\rm o}$ - $6^{\rm o}$).

When $B$ is a $G-C_0(X)$-algebra, where $X$ is a proper $G$-space, there exists also a canonical right action of $C_0(X)$ on $C^*(G,B)$ given by the right multiplication. Moreover, a similar right action of $C_0(X)$ exists also for any Hilbert module $\E$ over $C^*(G,B)$. It is given by $e\cdot f=\lim_{\alpha\to 0}e(e,e)((e,e)+\alpha)^{-1}f$, where $e\in \E,\,f\in C_0(X)$. In the next definition we will use this $C_0(X)$-action.

\smallskip

\begin{defi} Let $A$ and $B$ be $G-C_0(X)$-algebras. We define the group $KK^{G,X}(A,C^*(G,B))$ in the same way as the usual $KK^G(A,C^*(G,B))$, but with an additional assumption that for any triple $(\E,\psi:A\to\L(\E), F)$ and any $a\in A,\, f\in C_0(X)$, one has: $\psi(a)\cdot f=\psi(fa)$ in $\L(\E)$, where the right multiplication by $f$ on $\E$ is defined above. Note that when $G$ and $X$ are compact, $KK^{G,X}(C(X),C^*(G,B))$ is the same as $K_0^G(C^*(G,B))$.
\end{defi}
 
 \begin{prop}
 The group $KK^{G,X}(A,C^*(G,B))$ has the following additional properties. 

$1^{\rm o}$. For any $G-C_0(X)$-algebra $D$, there is a natural map: 
$$KK^{G,X}(A,C^*(G,B))\to KK^{G,X}(A\hotimes_{C_0(X)} D,C^*(G,B\hotimes_{C_0(X)} D))$$
given by the product with $1_D$.

\smallskip

$2^{\rm o}$. If $\alpha\in KK^{G,X}(A,C^*(G,B))$ and $\beta\in 
\cR KK^G(X;B,D)$, then the product $\alpha\otimes_{C^*(G,B)}j^G(\beta)$ belongs to $KK^{G,X}(A,C^*(G,D))$. 

\smallskip

$3^{\rm o}$. There is a natural product: 
$$KK^{G,X}(A,C^*(G,B))\otimes \cR KK^G(X;D,E)$$
$$\to KK^{G,X}(A\hotimes_{C_0(X)} D,C^*(G,B\hotimes_{C_0(X)} E)).$$
This product is commutative.

$4^{\rm o}$. There is an isomorphism:
$$KK^{G,X}(\Slf(X),C^*(G,C_0(X)))\simeq KK^{G,X}(C_0(X),C^*(G,\CG(X))).$$
\end{prop}

\Pf The first two statements are clear. The third one is a combination of the first two. For the commutativity of the product in $3^{\rm o}$, one has to identify $A\hotimes_{C_0(X)} D$ with $D\hotimes_{C_0(X)} A$ and $B\hotimes_{C_0(X)} E$ with $E\hotimes_{C_0(X)} B$ via $a_1\hotimes a_2\mapsto (-1)^{{\rm deg}a_1\cdot {\rm deg}a_2}a_2\hotimes a_1$. This is similar to \cite{Ka88}, theorem 2.14, $8^{\rm o}$. The last statement follows from theorem 7.6 and its proof. \cqfd

\subsection{Overview of definitions}

\emph{We will use the equivariant $KK^G$ in all notation for the symbol and index groups. In the non-equivariant (or stabilizer-invariant) case, one just has to drop the `$G$'.}

\medskip

In the definition of ellipticity we will follow \cite{Ka16}, section 3, for elliptic operators and \cite{Ka16}, section 6, for t-elliptic operators. As in \cite{Ka16}, we consider two options: 

\medskip

$1^{\rm o}$ Self-adjoint operators of degree $1$ on $\Z_2$-graded spaces have index in the $K^0$-groups.

\medskip

$2^{\rm o}$ Self-adjoint operators on ungraded spaces have index in the $K^1$-groups.

\begin{defi} {\bf (Elliptic case.)} A self-adjoint PDO acting on sections of a vector bundle $E$ is elliptic if its symbol $\sigma\in \L(C_0(p^*(E)))$ satisfies the condition: $f\cdot (\sigma^2-1)\in \K(C_0(p^*(E)))$ for any $f\in C_0(X)$ (both in the graded and non-graded case. Recall that $p:TX\to X$.)
\end{defi}

When the symbol $\sigma$ is elliptic, the coarse PDO construction of section 3 gives the index element in the group $K^*(C_0(X))$. If a locally compact group $G$ acts on $X$ properly and isometrically and the symbol $\sigma$ is elliptic and $G$-invariant then $[\F]=(L^2(E),\F)\in K^*_G(C_0(X))$.

Recall from \cite{Ka16} that the symbol group for elliptic operators is defined as $\cR KK^G_*(X;C_0(X),C_0(TX))$. If one wants to consider elliptic operators with coefficients in a unital $G$-algebra $B$ (see \cite{MiFo}), the symbol group will be $\cR KK^G_*(X;C_0(X),B\otimes C_0(TX))$ and the index group will be $KK^G_*(C_0(X),B)$. 

Note that in the case of a compact manifold $X$, the symbol group `shortens' to $K_*^G(C_0(TX))$. But it is important to emphasize that for a non-compact $X$, one can still consider symbols $\sigma$ belonging to $K_*^G(C_0(TX))$ if the group $G$ is compact. Such symbol has an additional property that the corresponding operator is Fredholm, because $\sigma^2-1$ is of strongly negative order (see 3.5 and 3.8, $4^{\rm o}$).

\begin{defi} {\bf (T-elliptic case.)} A self-adjoint PDO is t-elliptic if its symbol $\sigma\in \L(\gtS_{tr}(E))$ satisfies the condition: $f\cdot (\sigma^2-1)\in \K(\gtS_{tr}(E))$ for any $f\in C_0(X)$ (both in the graded and non-graded case).
\end{defi}

For t-elliptic operators (\cite{Ka16}, section 6), the symbol group is defined as $\cR KK^G_*(X;C_0(X),\gtS_{tr}(X))$, and in the case of operators with coefficients in a unital $G$-algebra $B$, the symbol group is $\cR KK_*^G(X;C_0(X),B\otimes \gtS_{tr}(X))$. Also for t-elliptic operators there are Clifford and tangent-Clifford symbol groups: $\cR KK_*^G(X;C_0(X),\CtG(X))$ and $\cR KK_*^G(X;C_0(X),\CG(TX))$ respectively. They are isomorphic to $\cR KK_*^G(X;C_0(X),\gtS_{tr}(X))$ by theorem 7.10.

For ``scalar'' t-elliptic operators (when $B=\C$), the index group was defined in \cite{Ka16}, proposition 6.4, as $K^*(C^*(G,C_0(X)))$. The corresponding index group for operators with coefficients in $B$ will be $KK_*^G(C^*(G,C_0(X)),B)$. 

The index of an operator $F$ is denoted $\ind(F)$. Another index $[F]\in K^*_G(\Slf(X))$, called analytical index, will be defined in 10.5. 

When the manifold $X$ and the group $G$ are compact, the symbol groups simplify to $K_*^G(\gtS_{tr}(X))\simeq K_*^G(\CtG(X))\simeq K_*^G(\CG(TX))$. However, when $G$ is compact but $X$ is not, we can still have symbols in the latter groups. The corresponding operators have a distributional index in $K^0(C^*(G))$ (see remark 10.4).

\begin{defi} {\bf (Leaf-wise elliptic case.)} A leaf-wise self-adjoint PDO is elliptic if its symbol $\sigma\in \L(\gtS_{lf}(E))$ satisfies the condition: $f\cdot (\sigma^2-1)\in \K(\gtS_{lf}(E))$ for any $f\in C_0(X)$ (both in the graded and non-graded case).
\end{defi}

For leaf-wise operators, the symbol group is $\cR KK_*^G(X;C_0(X),\gtS_{lf}(X))$, which is the same as
$\cR KK_*^G(X;C_0(X),C^*_\Gamma(X))$. The corresponding index group is $KK_*^{G,X}(C_0(X),C^*(G,C_0(X)))$.

When $G$ and $X$ are compact, the symbol group simplifies to $K_*^G(\gtS_{lf}(X))$, and for the index group one can take $K_*^G(C(X))$ (see remark 11.3). For compact $G$ and non-compact $X$, still symbols in $K_*^G(\gtS_{lf}(X))$ can be considered, with the index in $K_*^G(C_0(X))$ (remark 11.3).

\bigskip

\noindent{\bf \large Adding a new variable.}

\medskip

It is useful to generalize these definitions. Instead of just a `lonely' operator $\F$ acting on a Hilbert space or a Hilbert module we can consider the case when a $C_0(X)$-unital $G-C_0(X)$-algebra $A$ acts on the same Hilbert module and commutes modulo compact operators with $\F$.

\emph{A $C_0(X)$-algebra $B$ is called $C_0(X)$-\emph{unital} if there exists a subalgebra in $B$ isomorphic to $C_0(X)$, and the structure of the $C_0(X)$-algebra on $B$ is given by multiplication with elements of this subalgebra. (This means that all fibers $B_x,\, x\in X$, are unital.)}

\smallskip

For elliptic operators, the symbol group will be $\cR KK^G_*(X;A,B\otimes C_0(TX))$ and the index group $KK^G_*(A,B)$. 

\smallskip

For t-elliptic operators, we will get the symbol group $\cR KK^G_*(X;A,B\otimes \gtS_{tr}(X))$ and the index group $KK^G_*(C^*(G,A),B)$. 

\smallskip

For leaf-wise operators, the symbol group will be $\cR KK^G_*(X;A,\gtS_{lf}(X))$ and the index group $KK^{G,X}_*(A,C^*(G,C_0(X)))$.

\begin{ex} {\rm Let $B=\C$ and $A=\CtG(X)$. The Bott element $[\B_{\xi,\Gamma}]\in \cR KK^G(X;\CtG(X),\gtS_{tr}(X))$ defined before theorem 7.10 (see also definition 7.2 \cite{Ka16}) can be considered as an element of a t-elliptic symbol group. The corresponding element $[d_{X,\Gamma}]$ of the index group $K^0(C^*(G,\CtG(X)))$ is defined in \cite{Ka16}, 8.8. Here is this definition.

The Dirac element $[d_{X,\Gamma}]$ is given by the pair $(H,F_X)$, with $H=L^2(\Lambda^*(X))$, the action of $\CtG(X)$ on $H$ defined on (real) covectors by $\xi_1\oplus\xi_2\mapsto \ext(\xi_1)+\inter(\xi_1)+i(\ext(\xi_2)-\inter(\xi_2))$.
Here $\xi_1$ is a section of $\tau=T(X)$ and $\xi_2$ is a section of $\Gamma$. The representation of $C^*(G,\CtG(X))$ on $H$ is induced by the covariant representation defined by the action of $G$ and the above action of $\CtG(X)$, the operator $d_X$ is the exterior derivation operator, and $F_X=(d_X+d^*_X)(1+(d_X+d^*_X)^2)^{-1/2}$.}
\end{ex}

\begin{rem} {\rm In the t-elliptic case, the group $\cR KK_*^G(X; \gtS_{lf}(X), C_0(TX))$ can also be considered as the symbol group and $K^*_G(\gtS_{lf}(X))$ as the index group. See theorems 7.11 and 10.6.}
\end{rem}

\begin{ex} {\rm Another important element is the leaf-wise Dirac element $[D_{\Gamma,C_0(X)}]\in KK^{G,X}(\Slf(X),C^*(G,C_0(X)))$. Here is its definition:

Recall from the definition-lemma 7.3 the Hilbert $C_0(X)$-module $\T$ and the operator $\cC$ on it. The Hilbert module for our element will be $\Q=\T\hotimes_{C_0(X)} C^*(G,C_0(X))$. $\Slf(X)$ acts on it by multiplication over $\T$ as in 7.3. Let us present $\Q$ as the product $\{L^2(\tau_x)\}\hotimes_{C_0(X)} E\hotimes_{C_0(X)}C^*(G,C_0(X))$, where $E$ is the Hilbert module of continuous sections of $\Lambda^*(\tau)$. We consider the operator $\gtD_\Gamma$ as an operator on  $E\hotimes_{C_0(X)}C^*(G,C_0(X))$. The full operator for this element can be written as the $KK$-product $N_1(\cC\hotimes 1)+N_2(1\hotimes \tilde \gtD_\Gamma)$, where $\tilde \gtD_\Gamma$ is a $K$-theoretic connection for $\gtD_\Gamma$. Here $C_0(X)\cdot N_1\cdot \K(\T)\subset \K(\Q)$ and $C_0(X)\cdot N_2\cdot (1-\tilde \gtD_\Gamma^2)\subset \K(\Q)$. (The operator $1-\gtD_\Gamma^2$ has `negative order'.) Concerning the compactness of commutators between the operator $\cC$ and the Clifford variables of the operator $\gtD_\Gamma$, see the proof of 7.3.

From the element $[D_{\Gamma,C_0(X)}]$, we get the element $[D_{\Gamma,C_0(X)}]\otimes_{C_0(X)} 1_{\CG(X)}$ of the group $KK^{G,X}(\Slf(X)\hotimes_{C_0(X)}\CG(X),C^*(G,\CG(X)))$ (cf. proposition 8.2). By the $\cR KK$-duality theorem 7.6, this gives us the element $[D_\Gamma]\in KK^{G,X}(C_0(X),C^*(G,\CG(X)))$. The elements $[D_\Gamma]$ and $[D_{\Gamma,C_0(X)}]$ correspond to each other via the isomorphism of proposition 8.2, $4^{\rm o}$. The symbol corresponding to the element $D_\Gamma$ can be considered equal to the element $\gtf_\Gamma\in \cR KK^G(X;C_0(X),\CG(X)\otimes_{C_0(X)}\gtS_{lf}(X))$, although we cannot provide a \emph{direct} construction of the operator corresponding to this symbol.

The product of the element $[D_\Gamma]$ with $1_{\Ct(X)}$ in the sense of proposition 8.2, $1^{\rm o}$ will be denoted $[D_{\Gamma,\tau}]\in KK^{G,X}(\Ct(X),C^*(G,\CtG(X)))$.}
\end{ex}

\subsection{Theorems}

Now we come to a couple of important results. The first one is theorem 8.9 \cite{Ka16}:

\begin{theo} $[D_{\Gamma,\tau}]\otimes_{C^*(G,\CtG(X))}[d_{X,\Gamma}]=[d_X]$.
\end{theo}

\Pf The product in the statement is represented by the couple:
$$(\Ct(X)\to\L(\T\hotimes_{C_0(X)} L^2(\Lambda^*(X))), \cC \kkprod \,\gtf_\Gamma \kkprod \, \gtD_\Gamma \kkprod \, F_X).$$
The operator $\cC$, the $\xi$-multiplication part of the operator $\gtf_\Gamma$, and the Clifford variables of the operator $\gtD_\Gamma$ act on $\T$. The operator $F_X=(d_X+d^*_X)(1+\Delta)^{-1/2}$ (where $d_X$ is the exterior derivation, $\Delta$ is the Laplacian), as well as the Clifford variables of the operator $\gtf_\Gamma$, and the derivation variables of the operator $\gtD_\Gamma$, act on $L^2(\Lambda^*(X))$. The Clifford variables of the algebra $\Ct(X)$ also act on $L^2(\Lambda^*(X))$.

We first simplify this product a little bit by a rotation of Clifford variables. It is clear that the action of the Clifford variables of the operator $\gtD_\Gamma$ factors through the algebra $\M(\CG(X))$. The same is true for the operator $\gtf_\Gamma$. There exists such homomorphism of $Cl_{\Gamma\oplus\Gamma}(X)=\CG(X)\hotimes_{C_0(X)}\CG(X)$ into $\L(\T\hotimes_{C_0(X)} L^2(\Lambda^*(X))$ that the action of the Clifford variables of both these operators is contained in this action of $Cl_{\Gamma\oplus\Gamma}(X)$. By the rotation of $\Gamma_x\oplus\Gamma_x$ for all $x\in X$, we come to the situation when the Clifford variables of $\gtf_\Gamma$ act on $\T$ and the Clifford variables of $\gtD_\Gamma$ act on $L^2(\Lambda^*(X))$. This means that now the operators $\cC$ and $\gtf_\Gamma$ act on $\T$, and the operators $F_X$ and $\gtD_\Gamma$ act on $L^2(\Lambda^*(X))$. 

Now we can simplify both products $\cC\kkprod\gtf_\Gamma$ and $D_\Gamma\kkprod F_X$. The first one is homotopic to $\cC$ (see the proof of theorem 7.6). The second one is homotopic to $F_X$. The proof of this is given in the proof of theorem 8.9 \cite{Ka16}. It requires only a minor correction: one has to replace $\varphi_x(\xi)=f'_xf'^*_x(\xi)$ used in \cite{Ka16} with $\varphi_x(\xi)=(f'_xf'^*_x)^{1/2}(\xi)$ as it is defined at the beginning of section 5 of the present paper.

After all these changes we are left with the product $\cC\kkprod F_X$, and because $\cC$ represents the element $1_{C_0(X)}\in \cR KK(X;C_0(X),C_0(X))$, we can remove $\cC$. We arrive at the element $[d_X]$, as claimed. \cqfd

\medskip

The second result is part of theorem 8.15 \cite{Ka16}. Because of the confusion in \cite{Ka16} related with the $\cR KK$-duality of $\CG(X)$ with itself (instead of with $\Slf(X)$ as would be correct), we need to restate and reprove this important result.

\begin{theo} The following composition:
$$\cR KK^G_*(X; C_0(X),\CtG(X))\to KK^G_*(C^*(G,C_0(X)),C^*(G,\CtG(X)))$$
$$\to K_G^*(C^*(G,C_0(X)))\to K_G^*(\Slf(X))$$
coincides with the isomorphism of theorem 7.11.
Here the first map is $j^G$, the second one is the product with the element $[d_{X,\Gamma}]$, and the last one is the product with the element $[D_{\Gamma,C_0(X)}]$.
 \end{theo}

The proof is based on a lemma which copies lemma 8.16 \cite{Ka16}:

\begin{lem} Let $a\in \cR KK_*^G(X;C_0(X),\CtG(X))$. Then the triple product 
$$[D_{\Gamma,C_0(X)}]\otimes_{C^*(G,C_0(X))} j^G(a)\otimes_{C^*(G,\CtG(X))}[d_{X,\Gamma}]\in K^*_G(\Slf(X))$$ 
is equal to the product $\tilde a\otimes_{\Ct(X)}[d_X]$, where $\tilde a\in \cR KK^G_*(X;\Slf(X),\Ct(X))$ is $\cR KK$-dual to a.
\end{lem}

\noindent \emph{Proof of the lemma.} We will use the formula for $[d_X]$ from theorem 8.9. Then we need only to show that 
$$[D_{\Gamma,C_0(X)}]\otimes_{C^*(G,C_0(X))} j^G(a)=\tilde a \otimes_{\Ct(X)} [D_{\Gamma,\tau}].$$
The element $[D_{\Gamma,C_0(X)}]$ is the product of $1_{\gtS_{lf}(X)}\otimes_{C_0(X)}[D_\Gamma]$ and $j^G([d_{lf}])$, where $[d_{lf}]\in\cR KK(X;\CG(X)\otimes_{C_0(X)}\gtS_{lf}(X),C_0(X))$ is defined in 7.3 and used for $\cR KK$-duality in theorem 7.6. We have: $\tilde a=a\otimes_{\CG(X)} [d_{lf}]$ and $[d_{lf}]\otimes_{C_0(X)} a=1_{\CG(X)}\otimes_{C_0(X)} \tilde a$. Using this, it is easy to see that both sides of the  previous displayed formula are equal to the product of $[D_\Gamma]$ and $\tilde a$ in the sense of proposition 8.2, $3^{\rm o}$. \cqfd

\medskip

\noindent{\it Proof of theorem 8.10.} The composition in the statement of the theorem corresponds to the triple product from lemma 8.11. The isomorphism of theorem 7.11 corresponds to the product $\tilde a\otimes_{\Ct(X)}[d_X]$ from lemma 8.11. \cqfd

\section{Index theorem for elliptic operators}

The index theorem for elliptic operators is covered in great detail in section 4 of \cite{Ka16}, theorems 4.1-4.2. However, we suggest here another proof which will serve as a model for the proofs of index theorems for t-elliptic and leaf-wise elliptic operators in the next two sections.

Here is a reminder concerning the definition of the Dolbeault element (cf. \cite{Ka16}, 2.8). The manifold $TX$ is an almost complex manifold: the cotangent space $T^*_{(x,\xi)}(TX)$ at any point $(x,\xi)$ is the direct sum of the two orthogonal subspaces isomorphic to $T^*_x(X)$: the horizontal space $\tau_x^h=p^*(T^*_x(X))$ and the space dual to the vertical tangent space  $\tau_x^v=T_{\xi}(T_x(X))$. This allows to define the complex structure on the tangent space of the manifold $TX$ for any $(x,\xi)\in TX$. We will identify the complex cotangent bundle $T^*(TX)$ with the complexification of the lifted cotangent bundle of $X:\; p^*(T^*(X))\otimes \C$. The exterior algebra of this bundle is denoted by $\Lambda^*_\C(TX)$. The coordinates in the cotangent space $T^*_{(x,\xi)}(TX)$ will be denoted by $(\chi, \zeta)$. 

We define the symbol of the Dolbeault operator $\D$ as
$$\sigma_\D(x,\xi,\chi,\zeta)=\ext (-\chi+i\zeta)+\inter (-\chi+i\zeta)$$
on the exterior algebra bundle $\Lambda^*_\C(TX)$. (Note that $\inter=\ext^*$ in the Hermitian metric.) The differential operator $\D$ is essentially self-adjoint. The algebra $C_0(TX)$ acts on $L^2(\Lambda^*_\C(TX))$ by multiplication. The $K$-homology class $[\D_X]\in K^0_G(C_0(TX))$, called the Dolbeault element, is defined by the pair $(L^2(\Lambda^*_\C(TX)), \D(1+\D^2)^{-1/2})$.

\begin{theo} Let $X$ be a complete Riemannian manifold and $G$ a second countable locally compact group which acts on $X$ properly and isometrically. Let $F$ be a bounded, properly supported $G$-invariant elliptic operator on $X$ of order $0$ with the symbol $[\sigma_F]\in \cR KK^G_*(X;C_0(X),C_0(TX))$. Then the formula for the index $[F]$ of this operator is 
$$[F]=[\sigma_F]\otimes_{C_0(TX))} [\D_X]\in K^*_G(C_0(X)),$$
where $[\D_X]\in K^0_G(C_0(TX))$ is the Dolbeault element. 

If the group $G$ is compact and the symbol $[\sigma_F]$ belongs to $K^G_*(C_0(TX))$, then $F$ is Fredholm and defines an element of $R^*(G)=K_*^G(C({\rm point}))$ which is calculated by the same $KK$-product formula as above. 
\end{theo}

\Pf Both assertions of the theorem will be proved together. The fact that under the assumptions of the second assertion the operator $F$ is Fredholm follows from the coarse PDO construction (sections 2 and 3) because $||1-\sigma^2_F(x,\xi)||\to 0$ uniformly in $\xi$ when $x\to \infty$ (cf. corollary 2.6 and theorem 3.8).

On the manifold $TX$, let $\cC=\{\cC_x, x\in X\}$ be the family of Bott-Dirac operators on $\T=\{L^2(\Lambda^*(\tau_x))\}$ (see lemma 7.3), namely, $\cC_x=D_x(1+D_x^2)^{-1/2}$ where $D_x=(d_\xi+d^*_\xi+\ext(\xi)+\inter(\xi))$. Each of these operators has symbol 
$(\ext(\xi+i\zeta)+\inter(\xi+i\zeta))(1+||\xi||^2+||\zeta||^2)^{-1/2}$ 
and index $1$. Here $\zeta$ is the derivation variable, $\xi$ the multiplication variable. The family $\{\cC_x\}$ represents the element $1_{C_0(X)}\in KK^G(C_0(X),C_0(X))$. 
Taking the $KK$-product with $[F]\in K_G^*(C_0(X))$ 
we will get the same element $[F]\in K^*_G(C_0(X))$.

We want to prove that the $KK$-products 
$[\cC]\kkprod F$ and $\sigma_F\kkprod [\D_X]$ give the same element of $K_G^*(C_0(X))$. Both these products are represented on the Hilbert space $H=\T\otimes_{C_0(X)} L^2(E)$. Here we assume that $F$ acts on $L^2$-sections of the vector bundle $E$ over $X$. We can present the first product in the form: $(H, N_1(\cC\hotimes 1)+N_2(1\hotimes\tilde F))$, where $\tilde F$ is the $KK$-theoretic connection for $F$, and $N_1, N_2$ the operators entering the $KK$-product construction. According to lemma 3.11, we can take for $\tilde F$ a PDO on $TX$ with the symbol equal to the symbol of $F$ lifted to $TX$. 

The second product can be presented as $(H, N'_1\tilde\D_X+N'_2(1\hotimes \sigma_F))$, where $\tilde\D_X$ is the Dolbeault operator on $\{L^2(\Lambda^*(\tau_x))\}\otimes_{C_0(X)} L^2(X)$ lifted to $H$ (i.e. with the same symbol as the symbol of $\D_X$).

The homotopy that will join the two products will be the rotation homotopy which will go in the tangent spaces to $TX$ pointwise (see notation before the statement of the theorem). The rotation interchanging the derivation variable $\chi$ of $\tau_x^h$ and the multiplication variable $\xi$ of $\tau_x^v$ will be performed in $\tau_x^h\oplus \tau_x^v$.

We will present both products as pseudo-differential operators on $TX$ and produce a homotopy between their symbols. The symbol of $\cC \kkprod F$ on $TX$ can be written as 
$$N_1((\ext(\xi+i\zeta)+\inter(\xi+i\zeta))\hotimes 1)/(||\xi||^2+||\zeta||^2+1)^{1/2})+N_2(1\hotimes\sigma_F(x,\chi))$$
where 
$$N_1^2=(||\xi||^2+||\zeta||^2+1)/(||\xi||^2+||\chi||^2+||\zeta||^2+1),\;\;N_2^2=1-N_1^2.$$
The symbol of $\sigma_F\kkprod \D$ on $T^*(X)$ can be written as
$$N'_1((\ext(-\chi+i\zeta)+\inter(-\chi+i\zeta))\hotimes 1)/(||\chi||^2+||\zeta||^2+1)^{1/2}+N'_2(1\hotimes\sigma_F(x,\xi))$$
where 
$$N'^2_1=(||\chi||^2+||\zeta||^2+1)/(||\xi||^2+||\chi||^2+||\zeta||^2+1),\;\;N'^2_2=1-N'^2_1.$$ 

Here $\chi$ and $\zeta$ are the derivation variables, $x$ and $\xi$ the multiplication variables. 

The conditions required by definition 3.5 are easy to check. It is also easy to show that both symbols are elliptic in the sense that for each of these two symbols $\sigma$, on compact subsets of $X$ one has: $1-\sigma^2(x,\xi,\chi,\zeta)\to 0$ when $||\xi||^2+||\chi||^2+||\zeta||^2\to 0$. The $K$-theoretic connection property and the positivity property for the products (\cite{Ka88}, 2.10, (b) and (c)) are straightforward (cf. theorem 3.8).

The homotopy which joins the two operators with the symbols written above is given by the rotation in the $(\xi,\chi)$-variables: 
$$\xi\mapsto \cos t\cdot  \xi -\sin t\cdot \chi,\;\;\chi\mapsto \sin t\cdot \xi+ \cos t\cdot \chi,$$
with $0\le t\le \pi/2$. It is again easy to see that in the course of this homotopy the symbols remain elliptic. \cqfd

\section{Index theorems for t-elliptic operators}

Recall the definitions of the Clifford and the tangent-Clifford symbols for a t-elliptic operator $F$ with the symbol $[\sigma_F]\in \cR K_*^G(X;C_0(X),\gtS_{tr}(X))$(cf. \cite{Ka16}, definitions 8.11 and 8.13): 

\begin{defi} We define the tangent-Clifford symbol as
$$[\sigma_F^{tcl}]=[\sigma_F]\otimes_{\gtS_{tr}(X)}[\gtf_\Gamma]\in \cR KK^G_*(X;C_0(X), Cl_\Gamma(TX)),$$
where the element $\gtf_\Gamma\in \cR KK^G(X;\gtS_{tr}(X),\CG(TX))$ is defined before theorem 7.10. 

The Clifford symbol of a t-elliptic operator $F$ is defined as
$$[\sigma_F^{cl}]=[\sigma_F^{tcl}]\otimes_{\CG(TX)}[d_{\xi}]\in \cR KK^G_*(X;C_0(X), \CtG(X)),$$
where $[d_\xi]\in \cR KK^G(X;\CG(TX),\CtG(X))$ is also defined before theorem 7.10.
\end{defi}

\begin{rem}{\rm In the present paper we will not use the explicit formulas for the tangent-Clifford symbol given in \cite{Ka16}, 8.13. Instead, we will use the following formula for the tangent-Clifford symbol: 
$$\sigma_F^{tcl}(x,\xi)=\sigma_F(x,\xi)\hotimes 1+(1-\sigma_F^2)^{1/2}\hotimes \gtf_\Gamma(x,\xi)\in \L(C_0(p^*(E))\hotimes_{C_0(X)}\CG(X))$$
where $\gtf_\Gamma(x,\xi)=\varphi_x(\xi)(1+q_x(\xi))^{-1/2}$. Here $p:TX\to X$ is the projection and $E$ the vector bundle where the operator $F$ acts. We assume that $||\sigma_F||\le 1$ by cutting it above that level if necessary.   

Note that the exterior derivative $d_\xi (\sigma_F^{tcl}(x,\xi))$ vanishes at infinity in $\xi$ uniformly in $x\in X$ on compact subsets of $X$. Indeed, this is true for the $\sigma_F\hotimes 1$ summand because of the condition $2^{\rm o}$ of definition 7.1. Consequently, the same is true for $(1-\sigma_F^2)^{1/2}$ as well. Also it is easy to show that $||d_\xi(\gtf_\Gamma(x,\xi))||\le {\rm const}\cdot  (1+q_x(\xi))^{-1/2}$. In view of remark 7.2, this proves our claim. 

Note also that $1-(\sigma_F^{tcl})^2=(1-\sigma_F^2)(1+q_x(\xi))^{-1}$, and using remark 7.2, we see that $C_0(X)\cdot (1-(\sigma_F^{tcl})^2)$ is contained in $\K(C_0(p^*(E)))\hotimes_{C_0(X)}\CG(X)$. In fact, both $(1-\sigma_F^2)(1+q_x(\xi))^{-1}$ and $(1-\sigma_F^2)\varphi_x(\xi)(1+q_x(\xi))^{-1}$ vanish at infinity in $\xi$ uniformly on compact subsets in $x$ by 7.2. 

The only essential property of $\sigma_F^{tcl}(x,\xi)$ that is missing is the continuity in $x$ uniformly in $\xi$ (condition 3.5, $1^{\rm o}$).}
\end{rem}

The basic Dirac and Dolbeault elements for the next theorem are defined in \cite{Ka16}, 8.8 and 8.17. We have already recalled the definition of the Dirac element $[d_{X,\Gamma}]$ in example 8.6 above. Here is the definition of the Dolbeault elements:

The Clifford Dolbeault element $[\D^{cl}_{X,\Gamma}]\in K^0(C^*(G,\CG(TX)))$ (\cite{Ka16}, 8.17) comes from $[d_{X,\Gamma}]$ via the $\cR KK$-equivalence between $\CG(TX)$ and $\CtG(X)$ (see theorem 7.10). The explicit description of the element $[\D^{cl}_{X,\Gamma}]$ is as follows: Let $[\D_X]\in K^0_G(C_0(TX))$ be the Dolbeault element (see definition before theorem 9.1). We keep the Hilbert space and the operator $\D_X$ as in that definition, but extend the action of $C_0(TX)$ on $L^2(\Lambda^*_\C(TX))$ to the action of $Cl_\Gamma(TX)$. On (real) covectors of the Clifford part, this action is given by: $\chi\mapsto -(\ext (\chi)+\inter (\chi))$ (where $\chi$ is a section of $\Gamma \subset \tau=T(X)$). The resulting pair gives the element $[\D^{cl}_{X,\Gamma}]$. 

The Dolbeault element $[\D_{X,\Gamma}]\in K^0_G(C^*(G,\gtS_{tr}(X)))$ comes from $[\D^{cl}_{X,\Gamma}]$ via the $\cR KK$-equivalence between $\CG(TX)$ and $\gtS_{tr}(X)$ of theorem 7.10.

\smallskip

Recall from \cite{Ka16}, 8.1, that we consider $K^*(C^*(G,C_0(X)))$ as a subgroup of $K^*_G(C^*(G,C_0(X)))$ (see the beginning of subsection 8.1 for the definition of the $G$-action on $C^*(G,C_0(X))$). Namely, a triple $(\psi:C^*(G,C_0(X))\to\L(H),T)$ representing an element of $K^*(C^*(G,C_0(X)))$ (assuming that $\psi$ is non-degenerate) defines an element of $K_G^*(C^*(G,C_0(X)))$, where the $G$-action on $H$ is induced by $\psi$. 

\begin{theo} {\rm (Index theorem 8.18 \cite{Ka16}).} Let $X$ be a complete Riemannian manifold and $G$ a Lie group which acts on $X$ properly and isometrically. Let $F$ be a properly supported, $G$-invariant, $L^2$-bounded transversally elliptic operator on $X$ of order $0$. Then 
$$\ind(F)=j^G([\sigma_F^{cl}])\otimes_{C^*(G,\CtG(X))} [d_{X,\Gamma}]\in K^*(C^*(G,C_0(X))),$$ 
and
$$\ind(F)=j^G([\sigma_F])\otimes_{C^*(G,\gtS_{tr}(X))} [\D_{X,\Gamma}]\in K^*(C^*(G,C_0(X))),$$
and
$$\ind(F)=j^G([\sigma^{tcl}_F])\otimes_{C^*(G,Cl_\Gamma(TX))} [\D^{cl}_{X,\Gamma}]\in K^*(C^*(G,C_0(X))).$$
\end{theo}

\begin{rem}{\rm When the manifold $X$ and the group $G$ are compact, as noted already in subsection 8.2, the symbol group simplifies to $K_*^G(\gtS_{tr}(X))$, or to $K_*^G(\CtG(X))$, or to $K_*^G(\CG(TX))$. In this case, there exists the distributional index which is an element of $K^*(C^*(G))$. We will denote it $\ind^{\rm dist}(F)$. 

To obtain this index, one has to use the isomorphism   $K_*^G(\gtS_{tr}(X))\simeq KK_*^G(\C,\gtS_{tr}(X))$, apply the map $j^G$ and take the product with the element $\D_{X,\Gamma}$ as in the second formula of the theorem. Similarly can be modified the two other formulas of the theorem. This is a direct corollary of the theorem because the restriction map $KK_*^G(C(X),\gtS_{tr}(X)) \to KK_*^G(\C,\gtS_{tr}(X))$ for the inclusion $\C\subset C(X)$ corresponds to the restriction map $K^*(C^*(G,C_0(X)))\to K^*(C^*(G))$ for the inclusion $C^*(G)\subset C^*(G,C_0(X))$. 

In fact, it is enough to assume only that $G$ is compact, but not necessarily that $X$ is compact. Like in the case of elliptic operators (see theorem 9.1) if we assume that $[\sigma_F]$ is an element of $K_*^G(\gtS_{tr}(X))$, then $\sigma_F^2-1\in \gtS_{tr}(E)$. Following the coarse PDO construction of sections 2 and 3 (see in particular corollary 2.6 and theorem 3.8) and the proof of proposition 6.4 \cite{Ka16}, one can show that $\ind(F)\in K^*(C^*(G))$. This is the distributional index of $F$. The formulas of the theorem give the calculation of this distributional index as indicated above. The proof of this statement requires only minor modifications in the proof that follows.

We remark that in the case when $G$ is compact and $[\sigma_F]\in K_*^G(\gtS_{tr}(X))$, one can use the  homomorphism $K_*^G(\gtS_{tr}(X))\to KK^G_*(C(X^+),\gtS_{tr}(X))$, where $X^+$ is a one-point compactification of $X$. The formulas of the theorem applied to the symbol in the latter group give the index in $K^*(C^*(G,C(X^+)))$. This index, which we will denote $\Ind(F)$ allows to obtain both $\ind(F)$ and $\ind^{\rm dist}(F)$ by restricting $C(X^+)$ to $C_0(X)$ or, respectively, to the constant functions in $C(X^+)$.}
\end{rem}

\Pf All three formulas are equivalent. We will prove the third one. 

Let $F$ be a transversally elliptic operator of order $0$ acting on sections of the vector bundle $E$ over $X$. We will use the expression for the tangent-Clifford symbol given in remark 10.2. The product $j^G([\sigma_F^{tcl}])\otimes [\D^{cl}_{X,\Gamma}]$ is presented by the pair 
$$(C^*(G,C_0(X))\to \L(\E\hotimes_{\CG(X)} L^2(\Lambda^*(X))),N_1(\sigma_F^{tcl}\hotimes 1)+ N_2 (1\hotimes \tilde\D^{cl}_{X,\Gamma})),$$
where $\E$ is the Hilbert module of sections of the family of Hilbert spaces $\{E_x\hotimes Cl_{\Gamma_x}\otimes L^2(\tau_x)\}$ parametrized by  $x\in X$; $N_1, N_2$ are the operators involved in the $KK$-product construction, and $\tilde\D^{cl}_{X,\Gamma}$ is the $K$-theoretic connection for $\D^{cl}_{X,\Gamma}$. All Clifford variables act over $\Lambda^*(X)$.

We will follow closely the proof of theorem 9.1. In particular, we can choose the operators $N_1,N_2$ as in the proof of theorem 9.1. These will be the PDOs with the symbols:
$$\sigma_{N_1}=||\xi||(||\xi||^2+||\chi||^2+||\zeta||^2+1)^{-1/2}, \;\sigma_{N_2}=(1-\sigma_{N_1}^2)^{1/2}.$$
Here $\xi$ is the multiplication variable, $\chi$ and $\zeta$ are the derivation variables. The variable $\xi$ is used in $\sigma_F$ and $\gtf_\Gamma$. The variables $\chi$  and $\zeta$ are used in the symbol of $\D^{cl}_{X,\Gamma}$. 

One of the differences between the elliptic case (theorem 9.1) and the current t-elliptic one is that the symbol $\sigma_F^{tcl}$ is not a usual symbol of a PDO treated in section 3, although it has many similarities with such, as explained in remark 10.2. However, $\sigma_F^{tcl}$ is not continuous in $x$ uniformly in $\xi$, and the construction of the operator out of $\sigma_F^{tcl}$ requires the techniques of both sections 3 and 5. 

In order to deal with $\sigma_F^{tcl}$ along the rotation homotopy that will follow, we will treat the two summands of it differently. For example, $\sigma_{N_1}(\sigma_F\hotimes 1\hotimes 1)$ satisfies the conditions of definition 3.5, so this is a symbol of a PDO on the manifold $TX$. On the other hand, if the operator $(1-\sigma_F^2)^{1/2}\hotimes \gtf_\Gamma\hotimes 1$ commutes with some other PDO on a symbol level (like $N_1$ or $N_2$), then on the operator level the commutator becomes compact after multiplication by $C_0(X)$. This follows from proposition 3.3 (where we take $h=\varphi_x, k=(1-\sigma_F^2)^{1/2}$). The conditions of proposition 3.3 are satisfied - see remark 10.2.

The very specific feature of the t-elliptic case compared to the elliptic case of theorem 9.1 is that $\gtf_\Gamma(x,\xi)$ does not commute well with the operator $\tilde\D^{cl}_{X,\Gamma}$ because these two operators use the same Clifford variables. This deficiency of commutation is repaired by the presence of the algebra $C^*(G,C_0(X))$ which essentially provides an additional multiple $(1+\Delta_G)^{-1}$ killing the non-trivial commutators. (Here $\Delta_G$ is the orbital Laplacian, the operator with the symbol $q_x(\chi)$.)

We will use the rotation homotopy as in the proof theorem 9.1. The homotopy of all components of our operator, except $\gtf_\Gamma$, goes as in the proof of theorem 9.1. For the homotopy of $\gtf_\Gamma$, we use the rotation homotopy of the $\gtb$-symbol of the leaf-wise Dirac operator $\gtD_\Gamma$. The variable corresponding to $\xi$ for the $\gtb$-symbol will be $\eta$, as in definition 5.5. The analogous variable corresponding to $\chi$ will be called $\theta$. The rotation of the $\gtb$-symbol variables will go between $\eta$ and $\theta$.

\smallskip

Here is an important remark: 

Along this rotation homotopy, the commutator $[\gtf_\Gamma,\tilde\D^{cl}_{X,\Gamma}]$ will not be killed by multiplication with $C^*(G,C_0(X))$. The part of the commutator that is not killed is $\cos t \sin t [\varphi_x(\xi),\xi]$ multiplied by certain other positive operators. But since $[\varphi_x(\xi),\xi]\ge 0$, the commutator $[N_1 (\sigma_F^{tcl}\hotimes 1),N_2 (1\hotimes \tilde\D^{cl}_{X,\Gamma})]$ is also positive modulo compact operators and those operators which are killed by multiplication with $C^*(G,C_0(X))$. Therefore, in the process of the homotopy we get $(N_1 (\sigma_F^{tcl}\hotimes 1)+ N_2 (1\hotimes \tilde\D^{cl}_{X,\Gamma}))^2\ge 1$ modulo compact operators and those killed by multiplication with $C^*(G,C_0(X))$. So we can correct our homotopy by dividing $N_1 (\sigma_F^{tcl}\hotimes 1)+ N_2 (1\hotimes \tilde\D^{cl}_{X,\Gamma})$ by some positive invertible operator. At the end of the homotopy (at $t=\pi/2$) the non-trivial commutator again vanishes because $\sin t \cos t=0$.

\smallskip

After the rotation, $\sigma_F^{tcl}(x,\xi)$ becomes $\sigma_F^{tcl}(x,\chi)$. The operator $\sigma_F(x,\chi)$ is just $F$, and the operator $\gtf_\Gamma(x,\chi)$ corresponds to the leaf-wise Dirac operator $\gtD_\Gamma$ of subsection 5.3. Note that we assume that $||F||\le 1$ (along the whole homotopy) by cutting $F$ using functional calculus. Therefore $\sigma_F^{tcl}(x,\xi)$ turns into $T=F\hotimes 1+(1- F^2)^{1/2}\hotimes \gtD_\Gamma$ acting on $L^2(E\hotimes \Lambda^*(\tau))$.

The Dolbeault part is transformed by homotopy into the family of Bott-Dirac operators (of index $1$) $\cC=\{\cC_x, x\in X\}$ on $\{L^2(\Lambda^*(\tau_x))\}$, namely, $\cC_x=D_x(1+D_x^2)^{-1/2}$, where $D_x=(d_\xi+d^*_\xi+\ext(\xi)+\inter(\xi))$. We can remove the operator $\gtD_\Gamma$ from the operator $T$ because $\gtD_\Gamma$ is killed by multiplication with $C^*(G,C_0(X))$. After that we can remove the family $\cC$ as well because it corresponds to the element $1\in KK(C^*(G,C_0(X)),C^*(G,C_0(X))$. Therefore we arrive at the element $\ind(F)$. \cqfd

\medskip

The following definition replaces definition 8.10 \cite{Ka16} of the Clifford index of a t-elliptic operator.

\begin{defi} Let $F$ be a transversally elliptic operator with the index $\ind(F)\in K^*(C^*(G,C_0(X)))$, as defined in \cite{Ka16}, proposition 6.4. We will call $[F]=[D_{\Gamma,C_0(X)}]\otimes_{C^*(G,C_0(X))} \ind(F)
\in K^*_G(\gtS_{lf}(X))$ the analytical index of $F$. Here $[D_{\Gamma,C_0(X)}]\in KK^{G,X}(\Slf(X),C^*(G,C_0(X)))$ is the leaf-wise Dirac element (example 8.8).
\end{defi}

Recall that $\cR KK^G_*(X;C_0(X),\gtS_{tr}(X))\simeq \cR KK^G_*(X;\gtS_{lf}(X),C_0(TX))$, and both these groups are isomorphic to $K_G^*(\gtS_{lf}(X))$ (theorem 7.11). In particular, we can consider the symbol $[\sigma_F]$ of a t-elliptic operator $F$ acting on sections of a vector bundle $E$ over $X$ as an element $\cR KK^G_*(X;\gtS_{lf}(X),C_0(TX))$. 

\medskip 

The following theorem replaces the Clifford index theorems 8.12 and 8.14 of \cite{Ka16}.

\begin{theo} For a t-elliptic operator $F$ on a complete Riemannian manifold $X$ with a proper isometric action of the group $G$, the symbol $[\sigma_F]\in \cR KK^G_*(X;\gtS_{lf}(X),C_0(TX))$ and the analytical index $[F]\in K^*_G(\gtS_{lf}(X))$ are related by the isomorphism of theorem 7.11, namely, $[\sigma_F]\otimes_{C_0(TX)} [\D_{X}]=[F]$, where $[\D_{X}]\in K^0_G(C_0(TX))$ is the Dolbeault element.
\end{theo}

\Pf Combine theorems 8.10 and 10.3. \cqfd

\vfill
\eject

\noindent\emph{Sketch of the direct proof.}

\smallskip

We refer to the proof of theorem 10.3 for the details concerning the operator $T=F\hotimes 1+(1- F^2)^{1/2}(1\hotimes \gtD_\Gamma)$ on $L^2(E\hotimes \Lambda^*(\tau))$. Also recall from 7.3 the Hilbert $C_0(X)$-module $\T$ and the operator $\cC$. 

The analytical index $[F]$ is given by the pair $(\gtS_{lf}(X)\to \L(\T\otimes_{C_0(X)} L^2(E), A)$, where $A=N_1 (\cC\hotimes 1)+N_2(1\hotimes \tilde T)$. Here $N_1, N_2$ are the operators used in the $KK$-product construction, and $\tilde T$ is the $KK$-theoretic connection for $T$.  

The setting now is very similar to the one in the proof of theorem 10.3: the role of the algebra $C^*(G,C_0(X))$ is played here by $\gtS_{lf}(X)$. The operator $\gtD_\Gamma$ has Clifford variables which do not commute well with the operator $\cC$. This deficiency of commutation is repaired by the presence of the algebra $\gtS_{lf}(X)$ which provides an additional multiple killing the non-trivial commutators (see the proof of 7.3). 

We can perform the same rotation homotopy as in the proof of theorem 10.3 but going in the opposite direction, from $t=\pi/2$ to $t=0$. (The left action of the algebra $\gtS_{lf}(X)$ does not change along the homotopy.) In the process of this homotopy, there will be a non-trivial commutator between $\gtD_\Gamma$ and $\cC$, but its part related with $\xi$ will still be killed by multiplication with $\gtS_{lf}(X)$. The part of the commutator related with $\chi$ (the derivation part) will exist along the homotopy and it will be `positive', so the square of the whole product operator will be `$\ge 1$' (like in the proof of theorem 10.3). So one can `normalize' this homotopy. (See the proof of theorem 8.9 \cite{Ka16} for the positivity of the commutator.)

At the end of the homotopy, the family of operators $\cC$ will turn into the Dolbeault operator $\D_X$, and the operator $T$ will turn into the tangent-Clifford symbol $\sigma_F^{tcl}$. The $\varphi_x(\xi)(1+q_x(\xi))^{-1/2}$ part of $\sigma_F^{tcl}$ can be dropped because it is killed by multiplication with $\gtS_{lf}(X)$. Therefore we arrive at the product $[\sigma_F]\otimes_{C_0(TX)}[\D_X]$. \cqfd

\section{Index theorem for leaf-wise elliptic operators}

We will consider here the two cases: the one of $G$-invariant operators and the other of stabilizer-invariant operators (see assumption 5.4), simultaneously. We will keep the superscript `$G$' (which is needed in the first case and can be dropprd in the second one) in the notation of the $KK$-groups. 

\begin{theo} The symbol $[\sigma_F]\in \cR KK^G_*(X;C_0(X),\gtS_{lf}(X))$ and the index $\ind(F)\in KK^{G,X}_*(C_0(X),C^*(G,C_0(X)))$ of a leaf-wise elliptic operator $F$ are related by the formula: $\ind(F)=[\sigma_F]\otimes_{\gtS_{lf}(X)} [D_{\Gamma,C_0(X)}]$, where $[D_{\Gamma,C_0(X)}]\in KK^{G,X}(\gtS_{lf}(X),C^*(G,C_0(X)))$ is the leafwise Dirac element (see example 8.8 concerning $[D_{\Gamma,C_0(X)}]$ and proposition 8.2 concerning the $KK$-product).
\end{theo}

\Pf Recall that $\gtD_\Gamma$ is the leaf-wise operator on $\Lambda^*(\tau)$ with the symbol $\gtf_\Gamma(x,\xi)=\varphi_x(\xi)(1+q_x(\xi))^{-1/2}$ (see 5.12) . Also recall from example 8.8 the presentation $(\Slf(X)\to\L(\T\hotimes_{C_0(X)} C^*(G,C_0(X))), \cC \kkprod \gtD_\Gamma)$ of the element $[D_{\Gamma,C_0(X)}]$.

The product of $[\sigma_F]\in \cR KK(X;C_0(X),\gtS_{lf}(X))$ and $[D_{\Gamma,C_0(X)}]$ is given by the pair 
$$(\T\hotimes_{C_0(X)}E\hotimes_{C_0(X)}C^*(G,C_0(X)), N_1(\cC\hotimes 1\hotimes 1)+N_2(1\hotimes\sigma_F\hotimes 1)+N_3(1\hotimes 1\hotimes \tilde \gtD_\Gamma)).$$ 
Here $E$ is the Hilbert module where the operator $F$ acts, $\tilde \gtD_\Gamma$ is the $KK$-theoretic connection for $\gtD_\Gamma$, and $N_1,N_2,N_3$ are the operators used in the $KK$-product construction. 

Both $F$ and $\gtD_\Gamma$ admit a construction using the operator integration (subsection 5.4). We will apply the rotation homotopy to the operator representing the product. The symbols $\sigma_F$ and $\sigma_{\gtD_\Gamma}$ come from elements 
$$\gtb_F\in\L(E\hotimes_{C_0(X)}C_0(\gtG X))\;\; {\rm and} \;\;\gtb_{\gtD_\Gamma}\in\L(\Lambda^*(\tau)\hotimes_{C_0(X)} C_0(\gtG X))$$ via the map $\M(C_0(\gtG X))\to \M(\gtS_{lf}(X))$. Considering $\gtb_F\hotimes \gtb_{D_\Gamma}$ as an element of the tensor product 
$$\L(E\hotimes_{C_0(X)}\Lambda^*(\tau)\hotimes_{C_0(X)}C_0(\gtG X)\hotimes_{C_0(X)} C_0(\gtG X))$$
we will rotate (fiber-wise over $X$) the subspace $\gtG X$ of the fibered product $\gtG X\times_X \gtG X$ which is used for the operator integration. 

For any $x\in X$ and any $t,\,0\le t\le \pi/2$, we will choose in $\gtg_x\oplus \gtg_x$ the subspace isomorphic to $\gtg_x$ which corresponds to the second copy of $\gtg_x\oplus \gtg_x$ turned by the angle $t$, and we will apply the operator integration construction of subsection 5.4 using this turned second copy of $\gtg_x$. Note that there is no rotation related with the operator $\cC$. 

If we interpret the product operator written at the beginning as the one coming from the symbol $N'_1(\cC\hotimes 1\hotimes 1)+N'_2(1\hotimes\sigma_F\hotimes 1)+N'_3(1\hotimes 1\hotimes \gtf_\Gamma))$ by operator integration over the second copy of $\gtG X\times_X \gtG X$, then at the end of the homotopy we will get the operator coming from the same symbol by operator integration over the first copy of $\gtG X\times_X \gtG X$.

At the end of this homotopy (at $t=\pi/2$), we will get the operator $\tilde N_1(\cC\hotimes 1\hotimes 1)+\tilde N_2(1\hotimes F\hotimes 1)+\tilde N_3(1\hotimes 1\hotimes\gtf_\Gamma)$. As we have already seen in the proof of theorem 7.6, the product $\cC\kkprod\, \gtf_\Gamma$ is equal (homotopic) to $\cC$, so we can remove $[\gtf_\Gamma]$. Now we are left with the product (in the sense of proposition 8.2) of the element $[\cC]=1_{C_0(X)}$ and the element $\ind(F)$. So the result is $\ind(F)$. \cqfd

\begin{rem} {\rm If $X$ is a (complete Riemannian) $G$-manifold and $Y\subset X$ is a complete $G$-submanifold with boundary, we can restrict a leaf-wise operator $F$ from $X$ to $Y$ and the formula of theorem 11.1 will apply to $F|_Y$ (with $[\sigma_F]$ and $\ind(F)$ considered as elements of the corresponding groups on $Y$).}
\end{rem}

\begin{rem} {\rm When $G$ is compact, there is a natural map: $$KK^{G,X}_*(C_0(X),C^*(G,C_0(X)))\to K^G_*(C^*(G,C_0(X)))\to K_*^G(C_0(X)).$$ First arrow: for a compact $X$ -- see definition 8.1. The locally compact case follows by a one-point compactification. Second arrow: we can drop $G$ in $K^G$ and replace $K_*(C^*(G,C_0(X)))$ with $K_*^G(C_0(X))$. So in the case of a compact $G$, the index can be defined as an element of $K^G_*(C_0(X))$.}
\end{rem}

\begin{exs} {\rm $1^{\rm o}$ When $G$ and $X$ are compact and the action of $G$ on $X$ is free, we have:  $K_*^G(C(X))\simeq K_*(C(X/G))$. In this case, the index in $K_*(C(X/G))$ described above corresponds to the Atiyah-Singer index for a family of elliptic operators. So the family index theorem is a special case of the orbital index theorem 11.1.

$2^{\rm o}$ For each type of geometric Dirac operators listed in subsection 5.3 (see the `Construction' subtitle), theorem 11.1 gives a formula for the index. Certainly this index depends on the orbit structure of the manifold $X$. A simple concrete example of the index calculation will be given in subsection 12.3.}
\end{exs}

\section{Examples: t-elliptic and leaf-wise operators}

This section contains a discussion of some examples of t-elliptic and leaf-wise operators. T-elliptic operators are considered in subsections 12.1 and 12.2, a leaf-wise example is given in subsection 12.3. In all examples $G$ will be a compact group. 

Most t-elliptic operators that we consider will act on a Euclidean space $X=\R^{2n}$ or $X=\C^n$. In these cases, we will assume that $G$ acts on $X$ orthogonally, via a spin-representation (i.e via $G\to Spin(2n)\to SO(2n)$), or, in the case of $X=\C^n$, via a unitary representation. For t-elliptic operators in this section, we will be assuming that $[\sigma_F]\in K_0^G(\gtS_{tr}(X))$. 

\medskip

We will use the isomorphism: $K^0(C^*(G))=\widehat{R(G)}=\Hom (R(G),\Z)$. We consider $\Hom (R(G),\Z)$ as a module over $R(G)$, and the ring $R(G)$ as a submodule of $\widehat{R(G)}$. The $R(G)$-module structure on $K^*(C^*(G,B))$ for any $G$-algebra $B$ is given by the $KK$-product with the $j^G$-image of $R(G)$ in the group $KK(C^*(G,B),C^*(G,B))$.

\subsection{Atiyah's operators}

In his foundational lecture series \cite{Ati:TrEll}, M. Atiyah gave a very elaborate study of a class of t-elliptic operators on $\R^{2n}$. The index that he was calculating was the distributional index. We will indicate here a way of calculating the (topological) distributional index in certain examples based on theorem 10.3 and remark 10.4. 

\begin{ex} {\rm The most famous example of a t-elliptic operator for $G=S^1, X=\C$ is the differential operator $\tilde\partial_{\bar z}=\partial/\partial \bar z +z:L^2(\C^1)\to L^2(\C^1)\otimes\Lambda^1_\C(\C^1)$. We will consider a self-adjoint operator $\tilde\partial_{\bar z}+ \tilde\partial_{\bar z}^*$ with the symbol $\ext(z+i\xi)+\inter(z+i\xi)$ on $\Lambda_\C^*(\C^1)$ (where $\xi$ is a complex variable covector). We normalize it to order $0$ by dividing its symbol by $(1+||z||^2+||\xi||^2)^{-1/2}$. Let us denote the normalized operator by $F$. Its symbol $\sigma_F$ satisfies the strong decay condition at infinity of $X$: $||\sigma_F^2-1||_{z,\xi}\le {\rm const}\cdot (1+q_z(\xi))(1+||z||^2+||\xi||^2)^{-1}$, and defines an element $[\sigma_F]\in K_0^{S^1}(\gtS_{tr}(\C^1))$ (see \cite{Ka16}, 8.21).

Let us tensor $\Lambda_\C^*(\C^1)$ with another copy of $\Lambda^*_\C(\C^1)$ and consider the operator $(\tilde\partial_{\bar z}+ \tilde\partial_{\bar z}^*)\hotimes 1$ on $L^2(\C^1)\otimes\Lambda^*_\C(\C^2)$. There is a homotopy of this operator given by the homotopy of its symbol to an elliptic symbol: 
$$[\cos t (\ext z+\inter z)\hotimes 1+\sin t (1\hotimes (\ext z+\inter z))$$
$$+(\ext(i\xi)+\inter(i\xi))\hotimes 1](1+|z|^2+|\xi|^2)^{-1/2}$$ 
with $0\le t\le \pi/2$. At $t=\pi/2$, we get the symbol of the Bott-Dirac operator $\cC=[1\hotimes (\ext z+\inter z)+(\ext(i\xi)+\inter(i\xi))\hotimes 1](1+|z|^2+|\xi|^2)^{-1/2}$, which has index $1\in R(G)$.   

If we denote the basic representation of $G=S^1$ on $\C=\R^2$ by $\rho$, then the ($\Z_2$-graded) representation of $S^1$ on $\Lambda^*_\C(\C^1)$ corresponds to $[\Lambda^0]-[\Lambda^1]=1-\rho\in R(S^1)$. Note that the index of an elliptic operator can be calculated  by the formula of theorem 9.1, or theorem 10.3, and the results will be the same. Since we tensored our initial operator $F$ by $[\Lambda^*_\C(\C^1)]=1-\rho$, then using theorem 10.3 together with remark 10.4, we obtain: $1=\ind(\cC)=(1-\rho)\cdot \ind^{\rm dist}(F)\in R(S^1)$. 

Therefore, $\ind^{\rm dist}(F)=(1-\rho)^{-1}=\sum_{k=0}^\infty \rho^k \in \widehat {R(S^1)}$. This gives a calculation of the \emph{topological} distributional index for the operator $F$. (The calculation of the \emph{analytical} index of this operator is straightforward, see e.g. \cite{Ka16}, example 8.21. Of course, topological and analytical indices coincide.)} 
\end{ex}

\begin{ex} {\rm The previous example can be easily generalized to $X=\C^n$ with the action of $U(n)$. The operator will again be $\D_X+{\rm potential}$, where $\D_X$ is the Dolbeault operator and the potential is the operator $\ext(z)+\inter(z)$ with $z\in \Lambda^1(\C^n)=\C^n$. So the symbol of this operator is $\ext(z+i\xi)+\inter(z+i\xi)$ on $\Lambda_\C^*(\C^n)$ (where $\xi$ is a complex variable covector). We again normalize it to order $0$ by dividing its symbol by $(1+||z||^2+||\xi||^2)^{-1/2}$, and denote the normalized operator by $F$.

Tensoring $\Lambda_\C^*(\C^n)$ with another copy of $\Lambda_\C^*(\C^n)$ and doing the same homotopy of the operator as in example 12.1, we obtain the Bott-Dirac operator of index $1$. This time the element $[\Lambda_\C^*(\C^n)]=[\Lambda^{ev}]-[\Lambda^{od}]\in R(U(n))$. Its restriction to the maximal torus $T^n\subset U(n)$ is $\prod_{k=1}^n (1-\rho_k)$, which is invertible as an endomorphism of $\widehat {R(T^n)}$. The final answer is $\ind^{\rm dist}(F)=[\Lambda_\C^*(\C^n)]^{-1}$.}
\end{ex}

\subsection{Braverman's operators}

Let us consider now Braverman's operators (see \cite{Brav}).   The general setting is the following: $G$ is a compact Lie group, $X$ is a $G$-spin complete Riemannian manifold of dimension $2n$, $\D$ is a $G$-invariant Dirac operator acting on sections of the ($\Z_2$-graded) vector bundle $E$ (over $X$) which is a left module over $\Ct(X)$ (see, e.g., \cite{Ka16}, before definition 3.9). This is a differential operator with the symbol $\sigma_\D(x,\xi)=c(\xi)$, where $c$ means Clifford multiplication. When normalized to an operator $F$ of order $0$, it defines an element $[F]\in K^0_G(C_0(X))$. The symbol $[\sigma_F]$ of $F$ is an element of $\cR KK^G_0(X;C_0(X),C_0(TX))$, and the Clifford symbol $[\sigma^{cl}_F]$ is an element of $\cR KK^G_0(X;C_0(X),\Ct(X))$ (cf. \cite{Ka16}, 3.9, 3.10).

The Lie algebra of $G$ will be denoted $\gtg$. Braverman's `taming map' is an equivariant continuous map $\nu:X\to \gtg$ such that $||\nu(x)||\le 1$ and $||\nu(x)|| =1$ outside of a compact subset of $X$. It can be considered as a $G$-invariant section of the vector bundle $\gtg_X$ (see subsection 5.1). Using the $G$-action on $X$ (namely, the map $f':\gtg_X\to T(X)$), one gets a vector field $v(x)$ on $X$ (which will be used as a grading degree $1$ multiplier of $\Ct(X)$). This vector field has to be rescaled by multiplying it with a certain $G$-invariant positive scalar function $f(x)$ growing at infinity (see \cite{Brav}, definition 2.6). Then $\D_{fv}=\D+c(fv)$ is a Braverman type operator. Its symbol $\sigma_{\D_{fv}}(x,\xi)=c(\xi+f(x)v(x))$ is t-elliptic. When normalized to order $0$, it defines an element $[\sigma_{F_{fv}}]$ of the group $K_0^G(\Str(X))$. The corresponding index element will be denoted $\ind^{\rm dist}(F_{fv})$. Note that the vector field $v$ defines an element $[v]=(E,v)\in K_0^G(C_0(X))$. 

In the case of $X=\R^{2n}$, the Clifford module bundle $E$ is the product $X\times E_0$, where $E_0$ (see \cite{BGV}, 3.19) is the canonical Clifford module (isomorphic $\Lambda^*(\C^n)$ when the group $G$ acts via $U(n)$). The $G$-action on $E_0$ goes through the map $G\to Spin(2n)\to \End(E_0)$. Therefore, we have the corresponding element $[E_0^{ev}]-[E_0^{od}]\in R(G)$.

\begin{prop} Let $G=Spin(2n)$ and $D_{fv}$ a Braverman type operator on $\R^{2n}$. Then $([E_0^{ev}]-[E_0^{od}])\cdot [\sigma_{F_{fv}}]=[v]\otimes_{C_0(X)} [\sigma_F]\in K_0^G(\Str(X))$. Furthermore, $\ind^{\rm dist}(F_{fv})=([E_0^{ev}]-[E_0^{od}])^{-1}\cdot ([v]\otimes_{C_0(X)} [F])\in \widehat{R(G)}$.
\end{prop}

\noindent{\bf Remark.} The product $[v]\otimes_{C_0(X)} [\sigma_F]$ is defined by the following composition: $K_0^G(C_0(X))\otimes_{C_0(X)}\cR KK_0^G(X;C_0(X),C_0(TX))\to K_0^G(C_0(TX))\to K_0^G(\Str(X))$. The product $[v]\otimes_{C_0(X)} [F]\in R(G)$ is the usual pairing. 

\medskip

\Pf We tensor our operator with the space $E_0$ (as in examples 12.1 and 12.2). The resulting symbol $[\sigma_{F_{fv}}\hotimes 1]$ is homotopic to the product $[v]\,\kkprod \,[\sigma_F]$ by the homotopy similar to the one in examples 12.1 and 12.2. This gives the first formula. The second formula follows directly from this by the application of theorem 10.3 and remark 10.4. However, we need to invert the element $[E_0^{ev}]-[E_0^{od}]\in R(G)$. The maximal torus $T^n$ of $G=Spin(2n)$ acts on $E_0\simeq \Lambda^*(\C^n)$ in a coordinate-wise manner (as in example 12.2). Therefore, the restriction of $[E_0^{ev}]-[E_0^{od}]$ to the torus $T^n$ is $\prod_{k=1}^n (1-\rho_k)$, which is invertible as an endomorphism of $\widehat {R(T^n)}$.
\cqfd

\medskip

Proposition 12.3 resembles theorem 2.5 \cite{LRS} on the symbol level. Here is the exact analog of that theorem on the level of symbols for an arbitrary complete Riemannian manifold $X$. We consider the general setting described at the beginning of this subsection. Recall from \cite{LRS} that the vector field $v$ defines, in fact, an element $[v]\in K_0^G(\CG(X))$. The Clifford symbol $[\sigma_F^{cl}]$ of the (normalized) Dirac operator $F$ is an element of $\cR KK^G(X;C_0(X),\Ct(X))$ (see \cite{Ka16}, 3.8, 3.9, 3.10). It can be mapped into $\cR KK^G(X;\CG(X),\CtG(X))$ by tensoring it with $1_{\CG(X)}$. The Bott element $[\B_{\xi,\Gamma}]\in \cR KK^G(X;\CtG(X),\Str(X))$ is defined before theorem 7.10. 

\begin{prop} The symbol $[\sigma_{F_{fv}}]$ of the Braverman operator is equal to the triple product: $[v]\otimes_{C_0(X)}[\sigma_F^{cl}]\otimes_{C_0(X)}[\B_{\xi,\Gamma}]$.
\end{prop}

\Pf Let $E$ be the Clifford module bundle where the operator $\D$ acts. Using the calculation of $[\sigma_F^{cl}]$ in \cite{Ka16}, 3.9, 3.10, we easily see that the Hilbert module for the triple product is the image of the projection $P$ of \cite{Ka16}, 3.9, on $E\hotimes_{C_0(X)}\Lambda^*(X)$. But $\Lambda^*(X)\simeq E\hotimes_{C_0(X)}E$, so the image of $P$ is just $E$. The Bott operator $\beta$ acts on this $E$ as $c(\xi)(1+||\xi||^2)^{-1/2}$, the element $[v]$ acts as Clifford multiplication $c(v(x))$, so $[v]\kkprod [\beta]$ is obviously $[\sigma_{F_{fv}}]$. \cqfd

\subsection{A leaf-wise example}

Let the group $G=S^1$ act on the sphere $X=S^2$ by rotations around the axis which passes through the north and south poles. We will consider $S^2$ as $\C P^1$, and the (complex) vector bundle will be $\Lambda^*(\C P^1)$. The operator $D$ will be the leaf-wise Euler characteristic operator (called leaf-wise Dirac operator in subsection 5.3). However, unlike in subsection 5.3, here we will \emph{not} complexify $\Lambda^*(\C P^1)$, it is already a complex line bundle.

The north and south poles are fixed points. Outside of these two points, the action is free. One can write the exact sequence: 
$0\to C_0(X- 2\, {\rm pts})\to C(X)\to \C\oplus \C\to 0$. It is easy to check that the map $K_0^G(C(X))\to K_0^G(\C\oplus \C)$ is injective, and the quotient over the image is $\Z$. It is also well known that $K_0^{S^1}(\C P^1)=R(S^1)\oplus (1-\rho)R(S^1)$, and all maps in the exact sequence of the $K^G_*$-groups are $R(G)$-module maps.

The index is functorial. In each of the two poles, the index is $\pm(1-\rho)$. Therefore the index in $K_0^G(C(X))$ (see remark 11.3) is also $(\pm(1-\rho),\pm(1-\rho))$.

\appendix

\section{Appendix: Comments on my article \cite{Ka16}}

As stated in the Introduction, there were errors in my article \cite{Ka16}. In this Appendix, I will give an overview of noticed errors and indicate corrections.

\subsection{Section 5, \cite{Ka16}}

There are a couple of small technical errors in the proofs of theorem 5.8 and lemma 5.9.

\medskip

\noindent $1^{\rm o}$ \emph{Lemma 5.9}

\medskip

The error in the proof is the claim of the existence of a parametrix for the operator $T_{-\kappa,-\lambda}$ which has non-smooth coefficients. The first paragraph of the proof has to be replaced with the following:

\medskip

Let $B_\lambda$ be a properly supported $G$-invariant parametrix for the operator $D-i\lambda$, i.e. $(D-i\lambda)B_\lambda=1-S_\lambda$, where $S_\lambda$ is a properly supported smoothing operator. According to lemma 5.2, both $B_\lambda$ and $S_\lambda$ are bounded. We have: $T_{-\kappa,-\lambda}B_\lambda=1-S_\lambda-i\kappa\sigma_D(x,d\rho(x))B_\lambda$.

If a vector $u\in L^2(E)$ is orthogonal to the range of $T_{-\kappa,-\lambda}$, then for any $v\in C_c^\infty(E),\;(T_{-\kappa,-\lambda}(v),u)=0$. In particular, for any $v\in C_c^\infty(E)$, 
$$0=(T_{-\kappa,-\lambda}B_\lambda(v),u)=((1-S_\lambda-i\kappa\sigma_D(x,d\rho(x))B_\lambda)v,u)$$
$$=(v,(1-S^*_\lambda+i\kappa B^*_\lambda\sigma^*_D(x,d\rho(x))u),$$
which implies that $u=S^*_\lambda(u)-i\kappa\B^*_\lambda\sigma^*_D(x,d\rho(x))u$. Since $B^*_\lambda$ is a parametrix for $D+i\lambda$, it follows that $u$ is in the domain of $T_{\kappa,\lambda}$. From our initial assumption $(T_{-\kappa,-\lambda}(v),u)=0$, we deduce now that 
$(v,T_{\kappa,\lambda}u)=0$ for any $v\in C_c^\infty(E)$. so $u$ is in the kernel of $T_{\kappa,\lambda}$. 

\medskip

\noindent $2^{\rm o}$ \emph{Theorem 5.8}

\medskip

There was an erroneous choice of the function $ \mu(g)^{-1/2}\int_Xe^{-\kappa\rho(x)-\kappa\rho(g^{-1}x)}dx$ in the `End of the proof'. Here is the correction:

\medskip

\noindent \emph{End of the proof of theorem 5.8.} Let us choose $\kappa$ positive and large enough so that the following function of $g\in G:\; \mu(g)^{-1/2}\sup_Xe^{-\kappa\rho(x)-\kappa\rho(g^{-1}x)}$, belongs to $L^1(G)$. This is always possible. In fact, 
$e^{-\kappa\rho(x)-\kappa\rho(g^{-1}x)}\le e(g,x)=e^{-\kappa\rho(x,g^{-1}x)},$
where $\rho$ is the distance function. For a one orbit space $X$, $\mu(g)^{-1/2}e(g,x)$ obviously belongs to $L^1(G)$ if $\kappa$ is large enough. In general, under the assumptions that the $G$-action on $X$ is proper and isometric, and $X/G$ is compact, we can choose a compact subset $K\subset X$ such that $G\cdot K=X$, and the estimate is clear since $K$ has a finite diameter. Since $D$ is $G$-invariant and $X/G$ is compact, we can also choose $\lambda$ so that $2\kappa ||\sigma_D||<\lambda$. 

Let $v\in C_c^{\infty}(E)$. Because the operator $D\pm i\lambda=T_{0,\pm\lambda}$ is invertible on $L^2(E)$, one can solve the equation $(D\pm i\lambda)(u)=v$ in $L^2(E)$. But the operator $T_{\kappa,\pm\lambda}$ is also invertible on $L^2(E)$, and one can solve the equation $T_{\kappa,\pm\lambda}(u_1)=e^{\kappa\rho(x)}v$ as well. Calculating $(D \pm i\lambda)(u-e^{-\kappa\rho(x)}u_1)$ we get $0$. This means that $u=e^{-\kappa\rho(x)}u_1$, where $u_1\in L^2(E)$.

Let us calculate now $(u,u)_\E$. We have for any $g\in G$:
$$(u,u)_\E(g)=\mu(g)^{-1/2}\int_X(e^{-\kappa\rho(x)}u_1(x),e^{-\kappa\rho(g^{-1}x)}g(u_1)(x))dx$$
$$=\mu(g)^{-1/2}\int_Xe^{-\kappa\rho(x)-\kappa\rho(g^{-1}x)}(u_1(x),g(u_1)(x))dx$$
$$\le \mu(g)^{-1/2}\sup_X e^{-\kappa\rho(x)-\kappa\rho(g^{-1}x)} \cdot ||u_1||_{L^2(E)}^2.$$
The latter function of $g$ belongs to $L^1(G)$ by the previous choice of $\kappa$.

The same calculation shows that 
$$((1-a_\varepsilon) u,(1-a_\varepsilon)u)_\E(g)\le \mu(g)^{-1/2}\sup_Xe^{-\kappa\rho(x)-\kappa\rho(g^{-1}x)} \cdot ||(1-a_\varepsilon)u_1||_{L^2(E)}^2,$$
so when $\varepsilon\to 0$, the elements $a_{\varepsilon}u\in C_c(E)$ converge in $\E$, and the limit is $u$. But since $(D\pm i\lambda) (a_\varepsilon u)=a_{\varepsilon}(D\pm i\lambda)(u)-i\sigma_D(x,da_{\varepsilon})u$, we also get $(D\pm i\lambda) (a_\varepsilon u)\to (D\pm i\lambda) (u)=v$ (because $||\sigma_D(x,da_{\varepsilon})u||_\E\to 0$ by the same calculation as above, and $a_\varepsilon v=v$ for small $\varepsilon$ by the definition of $a_\varepsilon$). \cqfd

\subsection{Section 7, \cite{Ka16}}

\noindent $1^{\rm o}$ \emph{Definition 7.1}

\medskip

There is a major error at the end of definition 7.1. The claim was:

``There are two natural homomorphisms $\CG(X)\to \K(C_0(\Lambda^*_\Gamma(X)))$ given on continuous sections of $\Gamma$ by the maps $v\mapsto \ext(v)+\inter(v)$ and $v\mapsto i(\ext(v)-\inter(v))$. They generate the natural isomorphism: $\CG(X)\hotimes_{C_0(X)}\CG(X)\simeq \K(C_0(\Lambda^*_\Gamma(X))).$''

In fact, there is no such isomorphism in general unless $\Gamma$ is a vector bundle. This error led to a number of other errors in sections 7--9 of \cite{Ka16}. 

\medskip

\noindent $2^{\rm o}$ \emph{7.5--7.7}

\medskip

Essentially, the claim of the Remark 7.5 \cite{Ka16} was that $\CG(X)$ is $\cR KK$-dual to itself. This is wrong unless $\Gamma$ is a vector bundle. The correct statement is given in theorem 7.6 and corollary 7.7 of the present paper. Using the $\cR KK$-duality of the algebras $\CG(X)$ and $\Slf(X)$, one can create elements like $\Theta_{X,\Gamma}$ of definition 7.6 \cite{Ka16} and prove results like lemma 7.7 there. But we did not need this in the present paper.

\medskip

\noindent $3^{\rm o}$ \emph{7.8--7.9}

\medskip

These results are also wrong. They are replaced by the Poincare duality theorems 6.1 and 7.11 of the present paper.

\subsection{Section 8, \cite{Ka16}}

\noindent $1^{\rm o}$ \emph{Definition 8.3 and lemma 8.4}

\medskip

There are errors in definition 8.3 and lemma 8.4: that construction does not work.  The correct construction of the Dirac operator is given in subsection 5.3 of the present paper.

Note that the notation $\varphi_x(\xi)$ of the present paper is different from the same notation of \cite{Ka16}. The notation $\varphi_x(\xi)$ of the present paper corresponds to $\varphi_x^{1/2}(\xi)$ of \cite{Ka16} - see subsection 5.1 of the present paper.

\medskip

\noindent $2^{\rm o}$ \emph{Theorem 8.9}

\medskip

The proof of theorem 8.9 is essentially correct if one replaces $\varphi_x$ with $\varphi_x^{1/2}$. See details in the proof of theorem 8.9 of the present paper.

\medskip

\noindent $3^{\rm o}$ \emph{Definition 8.10}

\medskip

This definition is incorrect. It was replaced with the definition 10.5 of the present paper. 

\medskip

\noindent $4^{\rm o}$ \emph{Theorems 8.12 and 8.14}

\medskip

These theorems are incorrect. They are replaced with theorem 10.6 of the present paper. 

\medskip

\noindent $5^{\rm o}$ \emph{Theorem 8.15 and lemma 8.16}

\medskip

Theorem 8.15 is also incorrect as stated. The correct versions of theorem 8.15 and lemma 8.16 are theorem 8.10 and lemma 8.11 of the present paper.

\medskip

\noindent $6^{\rm o}$ \emph{Theorem 8.18}

\medskip

Theorem 8.18 is correct and is reproved as theorem 10.3 of the present paper.

\subsection{Section 9, \cite{Ka16}}

The results of this section need corrections very similar to the corrections indicated above for the section 8 \cite{Ka16}. Namely, definitions 9.2 and 9.3 should give the analytical index $[F]\in K^0(\gtS_{lf}(X))$ of a t-elliptic operator $F$. The precise implementation of this element is similar to the one given in the sketch of the direct proof of theorem 10.6 of the present paper. This leads to a corrected version of theorem 9.5 of \cite{Ka16}.


\begin{thebibliography}{99}

\bibitem{APT}
C. A. Akemann, G. K. Pedersen, J. Tomiyama:
\emph{Multipliers of \Cst-algebras.}
J. Funct. Analysis, 13 (1973), 277-301.

\bibitem{A-Sk II}
I. Androulidakis, G. Skandalis:
\emph{Pseudodifferential calculus on a singular foliation.}
J. Noncommut. Geom., 5 (2011), 125-152.

\bibitem{A-Sk III}
I. Androulidakis, G. Skandalis:
\emph{The analytic index of elliptic pseudodifferential operators on a singular foliation.}
J. K-Theory, 8 (2011), 363-385.

\bibitem{Ati:TrEll}
M. F. Atiyah: 
\emph{Elliptic operators and compact groups.}
Lect. Notes in Math., v. 401, Springer, 1974.

\bibitem{BGV}
N. Berline, E. Getzler, M. Vergne:
\emph{Heat kernels and Dirac operators.}
Springer, 1996.

\bibitem{Brav}
M. Braverman:
\emph{Index theorem for equivariant Dirac operators on non-compact manifolds.}
K-theory, 27 (2002), 61-101.

\bibitem{CM:L^2}
A. Connes, H. Moscovici: 
\emph{The $L^2$-index theorem for homogeneous spaces of Lie groups.} 
Annals of Math., 115 (1982), 291-330.

\bibitem{CS}
A. Connes, G. Skandalis:
\emph{The longitudinal index theorem for foliations.}
Publ. Res. Inst. Math. Sci. Kyoto Univ., 20 (1984), 1139-1183.

\bibitem{Gr80}
P. Green:
\emph{The structure of imprimitivity algebras.}
J. Funct. Anal., 36 (1980), 88-104.

\bibitem{GHT}
E. Guentner, N. Higson, J. Trout: 
\emph{Equivariant $E$-theory for $C^*$-algebras.} 
Mem. Amer. Math. Soc., 148 (2000).

\bibitem{Hor}
L. H\" ormander:
\emph{Fourier integral operators, I.}
Acta Math., 127 (1971), 79-183.

\bibitem{Hor_book}
L. H\" ormander:
\emph{The analysis of linear partial differential operators, III. Pseudo-differential operators.}
Springer-Verlag, 1985.

\bibitem{Ka75}
G. Kasparov:
\emph{Topological invariants of elliptic operators. I: $K$-homo\-logy.}
Izv. Akad. Nauk SSSR, Ser. Matem., 39 (1975), 796-838; translation: Math. USSR Izvestija, 9 (1975), 751-792.

\bibitem{KaJOT}
G. Kasparov:
\emph{Hilbert $C^*$-modules: Theorems of Stinespring and Voiculescu.}
J. Operator Theory, 4 (1980), 133-150.

\bibitem{Ka80}
G. Kasparov:
\emph{The operator $K$-functor and extensions of $C^{*}$-algebras.} 
Izvestiya Akad. Nauk SSSR, Ser. Matem., 44 (1980), 571-636; 
translation: Mathematics USSR - Izvesti\-ya, 16 (1981), 513-572.

\bibitem{Ka88} 
G. Kasparov: 
\emph{Equivariant $KK$-theory and the Novikov conjecture.} 
Invent. Math., 91 (1988), 147-201. 

\bibitem{Ka16}
G. Kasparov:
\emph{Elliptic and transversally elliptic index theory
from the viewpoint of $KK$-theory.}
J. Noncommut. Geom. 10 (2016), 1303-1378.

\bibitem{KS03}
G. Kasparov, G. Skandalis:
\emph{Groups acting properly on ``bolic'' spaces and the Novikov conjecture.}
Ann. of Math., 158 (2003), 165-206.

\bibitem{LRS}
Y. Loizides, R. Rodsphon, Y. Song:
\emph{A KK-theoretic perspective on deformed Dirac operators.}
Adv. Math., 380 (2021), 107604.

\bibitem{MiFo}
 A. S. Mishchenko, A. T. Fomenko:
 \emph{The index of elliptic operators over $C^*$-algebras.} Izv. Akad. Nauk SSSR Ser. Matem. 43 (1979), 831-859; translation: Mathematics USSR - Izvesti\-ya, 15 (1980), 87-112.
 
 \bibitem{Sk84}
G. Skandalis:
\emph{Some remarks on Kasparov theory.}
Journ. Funct. Anal., 56 (1984), 337-347.
 
 \bibitem{Sk88}
G. Skandalis:
\emph{Une notion de nucl\' earit\' e en $K$-th\' eorie.}
K-theory, 1 (1988), 549-573.
 
\bibitem{Thom_prepr}
K. Thomsen:
\emph{Asymptotic equivariant E-theory I.}
Preprint, 1997.


\end{thebibliography}
\end{document}